\documentclass{siamart220329}
\usepackage{amsmath,amsfonts,amssymb}

\usepackage{graphicx}                          %插图
\usepackage{setspace}  
\usepackage{caption}
\usepackage{listings}
\usepackage{tabularx}
\usepackage{xcolor}
\usepackage{indentfirst}
\usepackage{subcaption}

\usepackage{overpic}

\usepackage{verbatim}

\usepackage{tikz}

\usepackage{bm}

\usepackage[font=footnotesize,labelfont=bf]{caption}

\newtheorem{remark}{Remark}[section]

\newcommand\dd{\mathrm{d}}
\newcommand\pp{\partial}

\newcommand\e{\mathbf{e}}
\newcommand\x{\bm{x}}

\newcommand\uvec{\mathbf{u}}
\newcommand\X{\mathbf{X}}

\newcommand\y{\bm{y}}
\newcommand\z{\bm{z}}

\newcommand\n{{\bf n}}

\title{Energetic Variational Neural Network Discretizations of Gradient Flows}

\author{Ziqing Hu\thanks{Department of Applied and Computational Mathematics and Statistics, University of Notre Dame, 102G Crowley Hall, Notre Dame, IN 46556 USA (\email{zhu4@nd.edu}).} \and Chun Liu\thanks{Department of Applied Mathematics, Illinois Institute of Technology, Chicago, IL 60616, USA (\email{cliu124@iit.edu}).} \and Yiwei Wang \thanks{Corresponding author. Department of Mathematics, University of California, Riverside, Riverside, CA 92521, USA (\email{yiweiw@ucr.edu}).}
\and Zhiliang Xu\thanks{Department of Applied and Computational Mathematics and Statistics, University of Notre Dame, 102G Crowley Hall, Notre Dame, IN 46556 USA (\email{zxu2@nd.edu}).} }

\begin{document}

\maketitle

% \textcolor{red}{structure-preserving: variational structure. The scheme is energy-stable. computer memory efficient }

\begin{abstract}
  We present a structure-preserving Eulerian algorithm for solving $L^2$-gradient flows and a structure-preserving Lagrangian algorithm for solving generalized diffusions. Both algorithms employ neural networks as tools for spatial discretization.  Unlike most existing methods that construct numerical discretizations based on the strong or weak form of the underlying PDE, the proposed schemes are constructed based on the energy-dissipation law directly. This guarantees the monotonic decay of the system's free energy, which avoids unphysical states of solutions and is crucial for the long-term stability of numerical computations. To address challenges arising from nonlinear neural network discretization,  we perform temporal discretizations on these variational systems before spatial discretizations. This approach is computationally memory-efficient when implementing neural network-based algorithms. The proposed neural network-based schemes are mesh-free, allowing us to solve gradient flows in high dimensions. Various numerical experiments are presented to demonstrate the accuracy and energy stability of the proposed numerical schemes.
\end{abstract}

\section{Introduction}

Evolution equations with variational structures, often termed as gradient flows, have a wide range of applications in physics, material science, biology, and machine learning
\cite{doi2011onsager, peletier2014variational, weinan2020machine}.
These systems not only possess but also are determined by an energy-dissipation law, which consists of an energy of state variables that describes the energetic coupling and competition, and a dissipative mechanism that drives the system to an equilibrium. 

More precisely, beginning with a prescribed energy-dissipation law
\begin{equation} \label{ED_1}
    \frac{\dd}{\dd t} E^{\rm total} = - \triangle \leq 0~,
  \end{equation}
where $E^{\text{total}}$ is the sum of the Helmholtz free energy $\mathcal{F}$ and the kinetic energy $\mathcal{K}$, and $\triangle$ is the rate of energy dissipation, one could derive the dynamics, the corresponding partial differential equation (PDE), of the system by combining the Least Action Principle (LAP) and the Maximum Dissipation Principle (MDP). The LAP states that the equation of motion of a Hamiltonian system can be derived from the variation of the action functional $\mathcal{A} = \int_{0}^T (\mathcal{K} - \mathcal{F}) \dd t$ with respect to state variable $\x$. This gives rise to a unique procedure to derive the ``conservative force'' in the system. 
The MDP, on the other hand, derives the ``dissipative force'' by taking the variation of the dissipation potential $\mathcal{D}$, which equals to $\frac{1}{2}\triangle$ in the linear response regime (near equilibrium), with respect to $\x_t$. In turn, the force balance condition leads to the PDE of the system:
\begin{equation}\label{EnVarA_FB}
\frac{\delta \mathcal{D}}{\delta \x_t} = \frac{\delta \mathcal{A}}{\delta \x}~ .
\end{equation}
This procedure is known as the energetic variational approach (EnVarA) \cite{giga2017variational, liu2009energetic, wang2022entropy}, developed based on the celebrated work of Onsager \cite{onsager1931reciprocal,onsager1931reciprocal2} and Rayleigh \cite{strutt1871some} in non-equilibrium thermodynamics. During the past decades, the framework of EnVarA has shown to be a powerful tool for developing thermodynamically consistent models for various complex fluid systems, including two-phase flows \cite{xu2014energetic, X2017ModelPO, yue2004diffuse},  liquid crystals \cite{liu2009energetic}, ionic solutions \cite{eisenberg2010energy}, and reactive fluids \cite{wang2020field, wang2021two}.
In a certain sense, the energy-dissipation law (\ref{ED_1}), which comes from the first and second law of thermodynamics \cite{ericksen1998introduction}, provides a more intrinsic description of the underlying physics than the derived PDE (\ref{EnVarA_FB}), in particularly for systems that possess multiple structures and involve multiscale and multiphysics coupling. 

From a numerical standpoint, the energy-dissipation law (\ref{ED_1}) also serves as a valuable guideline for developing structure-preserving numerical schemes for these variational systems, as many straightforward PDE-based discretizations may fail to preserve the continuous variational structures, as well as the physical constraints, such as the conservation of mass or momentum, the positivity of mass density, and the dissipation of the energy.
As a recent development in this field \cite{liu2020lagrangian},
the authors proposed a numerical framework, called a \emph{discrete energetic variational approach}, to address these challenges. Similar methods are also used in \cite{QianXu02076A, xu2016variational}.
The key idea of this approach is to construct numerical discretizations directly based on the continuous energy-dissipation laws without using the underlying PDE (see Section 3.1 for details).
The approach has advantages in preserving a discrete counterpart of the continuous energy dissipation law, which is crucial for avoiding unphysical solutions and the long-term stability of numerical computations.  
Within the framework of the \emph{discrete energetic variational approach}, Eulerian \cite{noh2020dynamic}, Lagrangian \cite{liu2020lagrangian,liu2020variational}, and particle methods \cite{wang2021particle} have been developed for various problems. 

However, tackling high-dimensional variational problems remains a challenge. Traditional numerical discretization, such as finite difference, finite element, and finite volume methods, suffer from the well-known curse of dimensionality (CoD) \cite{Han_E_PNAS8505_18}. 
Namely, the computational complexity of the numerical algorithm increases exponentially as a function of the dimensionality of the problem \cite{Bellman_Dyn_Prog, Han_E_PNAS8505_18}. Although particle methods hold promise for addressing high-dimensional problems, standard particle methods often lack accuracy and are not suitable for problems involving derivatives.

During the recent years, neural networks (NNs) have demonstrated remarkable success across a wide spectrum of scientific disciplines \cite{devlin2018bert, hinton2012deep,  krizhevsky2012imagenet, lecun1998gradient,litjens2017survey}. Leveraging the potent expressive power of neural network \cite{ baron1993universal, cybenko1989approximation, E-CSIAM-AM-1-561}, particularly deep neural network (DNN) architectures, 
there exists a growing interest in developing neural network-based algorithms for PDEs, especially for those in high dimensions  \cite{dockhorn2019discussion, weinan2018deep, KHARAZMI2021113547,  khoo2021solving,  lagaris1998artificial,  raissi2017physics,  sirignano2018dgm, van1995neural,  wei2018machine}. 
Examples include  physics-informed neural network (PINN) \cite{raissi2017physics}, deep Ritz method (DRM) \cite{weinan2018deep},  deep Galerkin method (DGM) \cite{sirignano2018dgm}, variational PINN \cite{KHARAZMI2021113547}, and weak adversarial network (WAN) \cite{ZANG2020109409} to name a few. 
A key component of these aforementioned approaches is to represent solutions of PDEs via NNs. All of these approaches determine the optimal parameters of the NN by minimizing a loss function, which is often derived from either the strong or weak form of the PDE (see section 2.1 for more details).
By employing NN approximations, the approximate solution belongs to a space of nonlinear functions, which may lead to a robust estimation
by sparser representation and cheaper computation \cite{KHARAZMI2021113547}.

The goal of this paper is to combine neural network-based spatial discretization with the framework of the discrete energetic variational approach \cite{liu2020lagrangian}, to develop efficient and robust numerical schemes, termed as energetic variational neural network (EVNN) methods.  
To clarify the idea, we consider the following two types of gradient flows modeled by EnVarA:  
\begin{itemize}
\item {\bf $L^2$-gradient flow} that satisfies an energy-dissipation law
\begin{equation}\label{ED_L2GD}
\frac{\dd}{\dd t} \mathcal{F}[\varphi] =  - \int_{\Omega}  \eta(\varphi) |\varphi_t|^2 \dd \x~.
\end{equation}

\item {\bf Generalized diffusion} that satisfies an energy-dissipation law
\begin{equation}\label{ED_Diff}
   \frac{\dd}{\dd t} \mathcal{F}[\rho] = - \int_{\Omega} \eta(\rho) |\uvec|^2 \dd \x~,
\end{equation}
where $\rho$ satisfies $\rho_t + \nabla \cdot (\rho \uvec) = 0$~, known as the mass conservation.
\end{itemize}
Derivation of the underlying PDEs for these types of systems is described in  Sections \ref{sec:L2GD} and \ref{sec:gen_diff} of the paper. 

Our primary aim is to develop structure-preserving Eulerian algorithms to solve $L^2$-gradient flows and structure-preserving Lagrangian algorithms to solve generalized diffusions based on their energy-dissipation law by utilizing neural networks as a new tool of spatial discretization.
To overcome difficulties arising from neural network discretization, we develop a discretization approach that performs temporal discretizations before spatial discretizations. This approach leads to a computer-memory-efficient implementation of neural network-based algorithms.
Since neural networks are advantageous due to their ability to serve as parametric approximations for unknown functions even in high dimensions, 
the proposed neural-network-based algorithms are capable of solving gradient flows in high-dimension, such as these appear in machine learning \cite{EMaWu2020, huang2020convex, tabak2010density, wang2021particle}.

The rest of the paper is organized as follows. Section \ref{sec:prelim} reviews the EnVarA and some existing neural network-based numerical approaches for solving PDEs. Section \ref{sec:mthd} of the paper is devoted to the development of the proposed EVNN schemes for $L^2$-gradient flows and generalized diffusions. 
Various numerical experiments are presented in Section \ref{sec:simu} to demonstrate the
accuracy and energy stability of the proposed numerical methods. Conclusions are drawn in Section \ref{sec:con}.

\section{Preliminary}
\label{sec:prelim}

\subsection{Energetic Variational Approach}

In this subsection, we provide a detailed derivation of underlying PDEs for $L^2$-gradient flows and generalized diffusions by the general framework of EnVarA. We refer interested readers to \cite{giga2017variational, wang2022entropy} for a more comprehensive review of the EnVarA. In both systems, the kinetic energy $\mathcal{K} = 0$.

\subsubsection{$L^2$-gradient flow}\label{sec:L2GD}

$L^2$-gradient flows are those systems satisfying the energy-dissipation law:
\begin{equation}\label{GD_1}
\frac{\dd}{\dd t} \mathcal{F}[{\varphi}] = -  \int_{\Omega}  \eta  |\varphi_t|^2 \dd \x~.
\end{equation}  
where ${ \varphi}$ is the state variable, $\mathcal{F}[\varphi]$ is the Helmholtz free energy, and $ \eta > 0$ is the dissipation rate. % Gradient flows have natural variational structures in Eulerian coordinates. 
% We assume that $\varphi$ is a scalar function throughout this paper. 
One can also view ${ \varphi}$ as the generalized coordinates of the system \cite{doi2011onsager}. By treating $\varphi$ as $\x$, the variational procedure (\ref{EnVarA_FB}) leads to the following $L^2$-gradient flow equation:
\begin{equation*}
 \frac{\delta (\frac{1}{2} \int \eta |\varphi_t|^2 \dd \x )}{\delta \varphi_t}  = - \frac{\delta \mathcal{F}}{\delta \varphi}  \quad \Rightarrow \quad  \eta \varphi_t = - \frac{\delta \mathcal{F}}{\delta \varphi}~.
\end{equation*}
%Gradient flows have wide application in mathematical modeling. 
Many problems in soft matter physics, material science, and machine learning can be modeled as $L^2$-gradient flows. Examples include Allen--Cahn equation \cite{du2020phase}, Oseen--Frank and Landau--de Gennes models of liquid crystals \cite{de1993physics, lin2001static}, phase field crystal models of quasicrystal \cite{jiang2014numerical, lifshitz2007soft}, and generalized Ohta--Kawasaki model of diblock and triblock copolymers \cite{ohta1986equilibrium}.
For example, by taking $$ \mathcal{F}[\varphi] ={\displaystyle  \int_{\Omega} \frac{1}{2} |\nabla \varphi|^2 + \frac{1}{4 \epsilon^2} (\varphi^2 - 1)^2 \dd \x}~,$$ which is the classical Ginzburg--Landau free energy, and letting $\eta(\varphi) = 1$, 
one gets the Allen--Cahn equation $$ \varphi_t = \Delta \varphi - \frac{1}{\epsilon^2} (\varphi^2 - 1)\varphi\ .$$

\subsubsection{Generalized diffusion}
\label{sec:gen_diff}

Generalized diffusion describes the space-time evolution of a conserved quantity $\rho(\x, t)$.
Due to the physics law of mass conservation, $\rho(\x, t)$  satisfies the kinematics 
\begin{equation}\label{kinematics_1}
\pp_t \rho + \nabla \cdot (\rho \uvec) = 0\ ,
\end{equation} 
where $\uvec$ is an averaged velocity in the system.

\begin{figure}[!ht]
    \centering
    \includegraphics[width = 0.5 \linewidth]{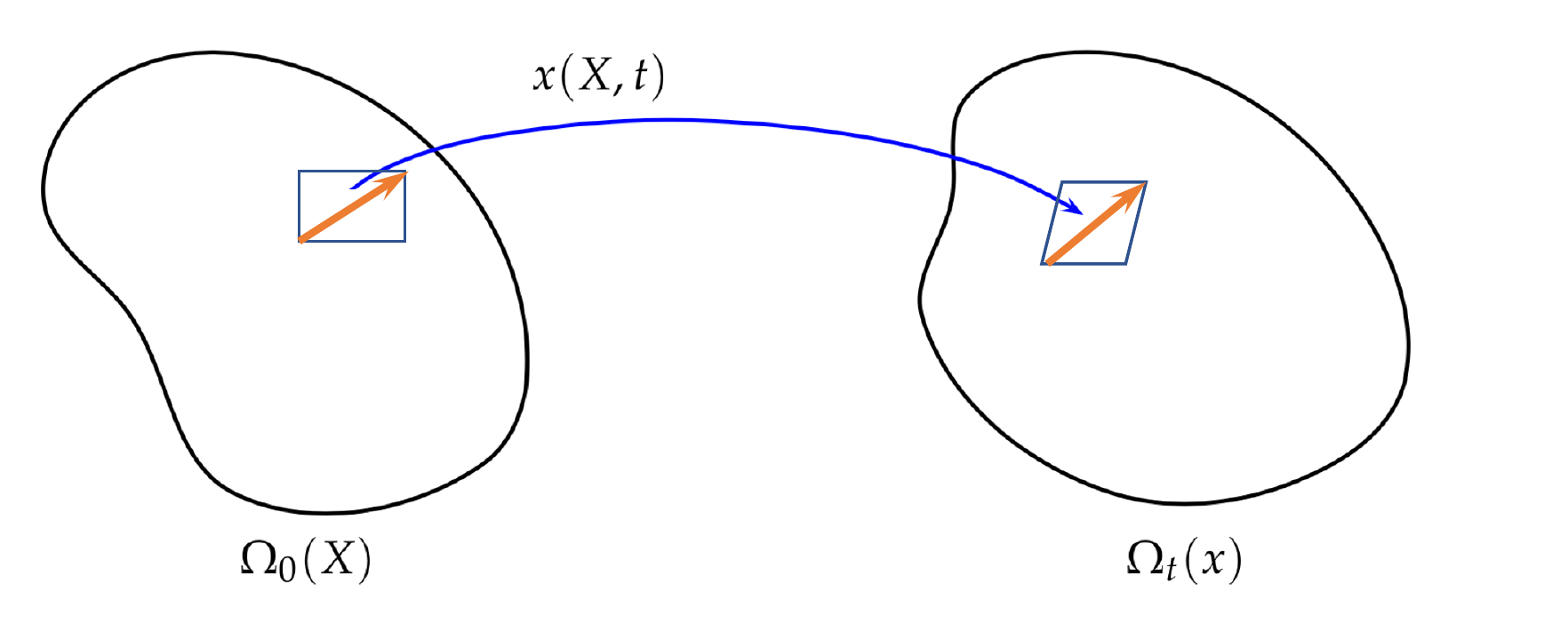}
    \caption{Schematic illustration of a flow map $x(\X, t)$.}
    \label{fig:flow_map}
\end{figure}

To derive the generalized diffusion equation by the EnVarA, one should introduce a Lagrangian description of the system. Given a velocity field $\uvec(\x, t)$, one can define a flow map $\x(\X, t)$ through
\begin{equation*}
\begin{cases}
    & \frac{\dd}{\dd t} \x(\X, t) = \uvec(\x(\X, t), t)~, \\
    & \x(\X, 0) = \X~, \\
\end{cases}
\end{equation*}
where $\X \in \Omega_0$ is the Lagrangian coordinates and $\x \in \Omega_t$ is the Eulerian coordinates. Here $\Omega_0$ is the initial configuration of the system, and $\Omega_t$ is the configuration of the system at time $t$. For a fixed $\X$, $\x(\X, t)$ describes the trajectory of a particle (or a material point) labeled by $\X$; while for a fixed $t$, $\x(\X, t)$ is a diffeomorphism from $\Omega_0$ to $\Omega_t$ (See Fig.~\ref{fig:flow_map} for the illustration.)  
The existence of the flow map $\x(\X, t)$ requires a certain regularity of $\uvec(\x,t)$, for instance, being  Lipschitz in $\x$. 

In a Lagrangian picture, the evolution of the density function $\rho(\x, t)$ is determined by the evolution of the $\x(\X, t)$ through the kinematics relation (\ref{kinematics_1}) (written in the Lagrangian frame of reference). 
More precisely, one can define the deformation tensor associated with the flow map $\x(\X, t)$ by
\begin{equation}
    {\sf\tilde{F}}(\x(\X,t),t) = {\sf F}(\X, t) = \nabla_{\X} \x(\X, t)~.
\end{equation}
Without ambiguity, we do not distinguish ${\sf F}$ and ${\sf \tilde{F}}$ in the rest of the paper. Let $\rho_0(\X)$ be the initial density, then the mass conservation indicates that
\begin{equation}
\rho(\x(\X, t), t) =  \rho_0(\X) / \det {\sf F}(\X, t)~, \quad \forall \X \in \Omega_0
\label{eq:gen_diff_Lag}
\end{equation}
in the Lagrangian frame of reference. This is equivalent to $\rho_t + \nabla \cdot(\rho \uvec) = 0$ in the Eulerian frame. 

Within (\ref{eq:gen_diff_Lag}), one can rewrite the energy-dissipation law (\ref{ED_Diff}) in terms of the flow map $\x(\X, t)$ and its velocity $\x_t(\X, t)$ in the reference domain $\Omega_0$.
A typical form of the free energy $\mathcal{F}[\rho]$ is given by
\begin{equation}\label{F_diffusion_general}
    \mathcal{F}[\rho] = \int_{\Omega} \omega(\rho) + \rho V(\x) +  \frac{1}{2}\left(\int_{\Omega} K(\x, \y) \rho(\x)\rho(\y) \dd \y \right) \dd \x~,
\end{equation}
where $V(\x)$ is a potential field, and $K(\x, \y)$ is a symmetric non-local interaction kernel. In Lagrangian coordinates, the free energy (\ref{F_diffusion_general}) and the dissipation in (\ref{ED_Diff}) become
\begin{equation*}
\begin{aligned}
    \mathcal{F}[\x(\X, t)] = &  \displaystyle  \int_{\Omega_0}  \omega\left(\tfrac{\rho_0}{\det {\sf F}}\right)  \det {\sf F} + V(\x) \rho_0(\X)  + \frac{\rho_0(\X)}{2} \left( \int_{\Omega_0} K(\x, \y) \rho_0(\X) \dd \X \right)   \dd \X \\
\end{aligned}
\end{equation*}
and $$\mathcal{D}(\x, \x_t) = \displaystyle \frac{1}{2} \int_{\Omega_0}  \eta(\rho_0 / \det {\sf F}) |\x_t|^2   \det {\sf F} \dd \X \ ,$$
By a direct computation, the force balance equation (\ref{EnVarA_FB}) can be written as (see \cite{giga2017variational} for a detailed computation)
% it is easy to compute that $\displaystyle \frac{\delta \mathcal{D}}{\delta \uvec} = \eta(\rho) \uvec$. So the force balance condition leads to 
\begin{equation}\label{FB1}
    \eta(\rho) \uvec = - \rho \nabla \mu, \quad  \mu =  \omega'(\rho) + V(\x) + K * \rho.
\end{equation}
In Eulerian coordinates, combining (\ref{FB1}) with the kinematics equation (\ref{kinematics_1}), one obtains a generalized diffusion equation
\begin{equation}
\label{G_diffusion}
\rho_t = \nabla \cdot \left(  m(\rho)  \nabla \mu\right)~, \quad m(\rho) = \rho^2 / \eta(\rho) \ ,
\end{equation}
where $m(\rho)$ is known as the mobility in physics. Many PDEs in a wide range of applications, including the porous medium equations (PME) \cite{vazquez2007porous}, nonlinear Fokker--Planck equations \cite{jordan1998variational}, Cahn--Hilliard equations \cite{liu2020lagrangian}, Keller--Segel equations \cite{keller1970initiation}, and Poisson--Nernst--Planck (PNP) equations \cite{eisenberg2010energy}, can be obtained by choosing $\mathcal{F}[\rho]$ and $\eta(\rho)$ differently.

The velocity equation (\ref{FB1}) can be viewed as the equation of the flow map $\x(\X, t)$ in Lagrangian coordinates, i.e.,
\begin{equation}\label{Flow_eq}
  \eta \left( \frac{\rho_0(\X)}{ \det {\sf F}(\X)}  \right)   \frac{\dd}{\dd t} \x(\X, t) = - \frac{\rho_0(\X)}{\det {\sf F}(\X) }  {\sf F}^{-\rm T}\nabla_{\X} \mu \left(\frac{\rho_0(\X)}{\det {\sf F}(\X)} \right)~,
\end{equation}
where ${\sf F}(\X) = \nabla_{\X} \x(\X, t)$. The flow map equation (\ref{Flow_eq}) is a highly nonlinear equation of $\x(\X, t)$ that involves both ${\sf F}$ and $\det {\sf F}$. It is rather unclear how to introduce a suitable spatial discretization to Eq.~(\ref{Flow_eq}) by looking at the equation.
However, a generalized diffusion can be viewed as an $L^2$-gradient flow of the flow map in the space of diffeomorphism. This perspective gives rise to a natural discretization of generalized diffusions \cite{liu2020lagrangian}. 

% In the rest of the paper, we always take $\eta(\rho) = \rho$ for simplicity if not stated differently.

\begin{remark}
In the case that $\eta(\rho) = \rho$, a generalized diffusion can be viewed as a Wasserstein gradient flows in the space of all probability densities having finite second moments $\mathcal{P}_2(\Omega)$ {\rm \cite{jordan1998variational}}. Formally, the Wasserstein gradient flow can be defined as a continuous time limit ($\tau \rightarrow 0$) from a semi-discrete JKO scheme,
\begin{equation}\label{JKO}
  \rho^{n+1} = \mathop{\arg\min}_{\rho \in \mathcal{P}_2(\Omega)} \frac{1}{2 \tau} W_2 (\rho, \rho^n)^2 + \mathcal{F}[\rho]~, \quad n = 0, 1, 2\ldots~,
\end{equation}
where $\displaystyle \mathcal{P}_2 (\Omega) = \left\{ \rho: \Omega \rightarrow [0, \infty) ~|~ \int_{\Omega} \rho \ \dd \x = 1, ~~\int_{\Omega} |\x|^2 \rho(\x) \dd \x < \infty \right\}$ and $W_2(\rho,\rho^n)$ is the Wasserstein distance between $\rho$ and $\rho^n$. 
The Wasserstein gradient flow is an Eulerian description of these systems {\rm \cite{benamou2016augmented}}.
Other choices of dissipation can define other metrics in the space of probability measures {\rm \cite{adams2013large, lisini2012cahn}}. 

\end{remark}

%\subsection{Eulerian, Lagrangian, and Particle methods}
% \textcolor{red}{move to introduction}

\subsection{Neural-network-based numerical schemes for PDEs}

In this subsection, we briefly review some existing neural network-based algorithms for solving PDEs. We refer interested readers to \cite{weinan2021algorithms,lu2021deepxde} for detailed reviews.

Considering a PDE subject to a certain boundary condition
\begin{equation}\label{PDE}
  \begin{aligned}
& \mathcal{L} \varphi(\x) = f(\x)~, ~~ \x \in \Omega \subset \mathbb{R}^d ~, \quad  \mathcal{B} \varphi(\x) = g(\x)~, ~~ \x \in \pp \Omega  ~,
  \end{aligned}
\end{equation}
and assuming the solution can be approximated by a neural network $\varphi_{\rm NN}(\x; {\bm \Theta})$, where ${\bm \Theta}$ is the set of parameters in the neural network,
the Physics-Informed Neural Network (PINN) \cite{raissi2017physics} and similar methods \cite{dockhorn2019discussion, sirignano2018dgm, wei2018machine} find the optimal parameter ${\bm \Theta}^*$ by minimizing a loss function defined as:
\begin{equation}
L({\bm \Theta}) = \frac{1}{N_{in}} \sum_{i = 1}^{N_{in}} \left( \mathcal{L}\varphi_{\rm NN} (\x_i, {\bm \Theta}) - f(\x_i) \right)^2 + \frac{\lambda}{N_b} \sum_{j = 1}^{N_b} (\mathcal{B} \varphi_{\rm NN}(\bm{s}_j; {\bm \Theta}) - g(\bm{s}_j))^2 ~.
\label{eq:strong_form_loss}
\end{equation}
Here $\{\x_i\}_{i=1}^{N_{in}}$ and $\{ \bm{s}_j \}_{j=1}^{N_b}$ are sets of samples in $\Omega$ and $\pp \Omega$ respectively, which can be drawn uniformly or by following some prescribed distributions.  The parameter $\lambda$ is used to weight the sum in the loss function.
Minimizers of the loss function (\ref{eq:strong_form_loss}) can be obtained by using some optimization methods, such as AdaGrad and Adam.
Of course, the objective function of this minimization problem in general is not convex even when the initial
problem is. Obtaining the global minimum of (\ref{eq:strong_form_loss}) is highly non-trivial.

In contrast to the PINN, which is based on the strong form of a PDE, the deep Ritz method (DRM) \cite{weinan2018deep}
is designed for solving PDEs using their variational formulations. % 
If Eq.~(\ref{PDE}) is an Euler-Lagrangian equation of some energy functional 
\begin{equation}\label{DRM_energy}
I[\varphi, \nabla \varphi] =  \int_{\Omega} W(\varphi(\x), \nabla_{\x} \varphi (\x)) \dd \x + \lambda \int_{\pp \Omega} |\mathcal{B} \varphi - g(\x) |^2   \dd S~, 
\end{equation}
then the loss function can be defined directly by
\begin{equation}
L(\theta) = \frac{1}{N_{in}} \sum_{i = 1}^{N_{in}}  W(\varphi_{\rm NN}(\x_i), \nabla_{\x} \varphi_{\rm NN}(\x_i)) + \frac{\lambda}{N_b} \sum_{j = 1}^{N_b} (\mathcal{B}\varphi_{\rm NN}(\bm{s}_j) - g(\bm{s}_j))^2~.
\end{equation}
The last term in (\ref{DRM_energy}) is a penalty term for the boundary condition. The original Dirichlet boundary condition can be recovered with $\lambda \rightarrow \infty$.
 Again, the samples $\{\x_i\}_{i=1}^{N_{in}}$ and $\{ \bm{s}_j \}_{j=1}^{N_b}$ can be uniformly sampled from $\Omega$ and $\pp \Omega$ or sampled following some other prescribed
distributions, respectively. Additionally, both the variational PINN \cite{KHARAZMI2021113547} and the WAN \cite{ZANG2020109409} utilized the Galerkin formulation to solve PDEs. The variational PINN  \cite{KHARAZMI2021113547}, stems from the Petrov--Galerkin method, represents the solution 
via a DNN, and keeps test functions belonging to linear function spaces. The WAN employs the primal and adversarial networks to parameterize the weak solution and test functions respectively  \cite{ZANG2020109409}. 

%This is also the starting point of the ``Deep BSDE method'' developed in \cite{weinanE_Han_Jentzen17, Han_E_PNAS8505_18}  to solve high-dimensional nonlinear parabolic PDEs.
%Results from  %\cite{weinanE_Han_Jentzen17, Han_E_PNAS8505_18} also 
%these works underscore the importance of formulating PDEs as variational problems to achieve better performance.

The neural network-based algorithms mentioned above focus on elliptic equations. For an evolution equation of the form:
\begin{equation}
\left\{
\begin{aligned}
& \pp_t \varphi(\x, t) =  F(t, \x, \varphi)~, \quad (\x, t) \in \Omega \times (0, \infty)~, \\
& \varphi(\x, 0) = \varphi_0(\x)~,\quad \quad \quad \x \in \Omega \\
\end{aligned}
\right.
\end{equation}
with a suitable boundary condition, the predominant approach is to treat the time variable as an additional dimension, and the PINN type techniques can be applied.
However, these approaches are often expensive and may fail to capture the dynamics of these systems. Unlike spatial variables, there exists an inherent order in time, where the solution at time $T$ is determined by the solutions in $t < T$.
It is difficult to generate samples in time to maintain this causality.
Very recently, new efforts have been made in using NNs to solve evolution PDEs \cite{bruna2022neural, du2021evolutional}.
The idea of these approaches is to use NNs with time-dependent parameters to represent the solutions. For instance, Neural Galerkin method, proposed in \cite{bruna2022neural}, parameterizes the solution as $\varphi_h(\x; {\bm \Theta}(t))$ with ${\bm \Theta}(t)$ being the NN parameters, and defines a loss function in terms of ${\bm \Theta}$ and $\dot{\bm \Theta}$ through a residual function, i.e.,
\begin{equation}
J({\bm \Theta}, {\bm \eta}) = \frac{1}{2} \int | \nabla_{\bm \Theta} \varphi_h  \cdot {\bm \eta} - F(\x, t, \varphi_h(\x; {\bm \Theta}))    |^2   \dd \nu_{\Theta} (\x)~,
\end{equation}
where $\nu_{\Theta} (\x)$ is a suitable measure which might depend on ${\bm \Theta}$.
By taking variation of $J({\bm \Theta}, {\bm \eta})$ with respect to ${\bm \eta}$, the neural Galerkin method arrives at an ODE of ${\bm \Theta}(t)$:
\begin{equation}\label{ODE_NGM}
   {\sf  M}( {\bm \Theta} ) \dot{\bm \Theta} = F(t, {\bm \Theta})~, \quad {\bm \Theta}(0) = {\bm \Theta}_0~,
\end{equation}
where ${\sf M}( {\bm \Theta}) = \int  \nabla_{\bm \Theta} \varphi_h \otimes  \nabla_{\bm \Theta} \varphi_h \dd \nu_{\bm \Theta} (\x),$ and  $F(t, {\bm \Theta}) = \int   \nabla_{\bm \Theta} \varphi_h F(\x, t, \varphi_h )  \dd \nu_{\bm \Theta} (\x).$ The ODE (\ref{ODE_NGM}) can be solved by standard explicit or implicit time-marching schemes.

It's important to note that, with the exception of DRM, all existing neural-network-based methods are developed based on either the strong or weak forms of PDEs. While DRM utilizes the free energy functional, it is only suitable for solving static problems, i.e., finding the equilibrium of the system. Additionally, most of these existing approaches are Eulerian methods. For certain types of PDEs, like generalized diffusions, these methods may fail to preserve physical constraints, such as positivity and the conservation of mass of a probability function.

\section{Energetic Variational Neural Network}
\label{sec:mthd}

In this section, we present the structure-preserving EVNN discretization for solving both $L^2$-gradient flows and generalized diffusions, As mentioned earlier, the goal is to construct a neural-network discretization based on the energy-dissipation law, without working on the underlying PDE. One can view our method as a generalization of the DRM to evolution equations.

\subsection{EVNN scheme for $L^2$-gradient flow}

Before we discuss the neural network discretization, we first briefly review the discrete energetic variational approach proposed in \cite{liu2020lagrangian}. Given a continuous energy-dissipation law (\ref{ED_1}), the discrete energetic variational approach first constructs a finite-dimensional approximation to a continuous energy-dissipation law by introducing a spatial discretization of the state variable $\varphi$, denoted by $\varphi_h(\x; {\bm \Theta}(t) )$, where ${\bm \Theta}(t) \in \mathbb{R}^K$ is the parameter to be determined. By replacing $\varphi$ by $\varphi_h$, one obtain a semi-discrete energy-dissipation law (here we assume $\mathcal{K}$ in the continuous model for simplicity) in terms of ${\bm \Theta}$:
% \label{ED_L2GD}
\begin{equation}
\label{discrte_En_1}
\frac{\dd}{\dd t} \mathcal{F}_h \left[ {\bm \Theta}(t) \right] = -  \triangle_h [{\bm \Theta}(t), {\bm \Theta}'(t)]~,
\end{equation}
where $\mathcal{F}_h \left[ {\bm \Theta}(t) \right] = \mathcal{F}[\varphi_h(\x;  {\bm \Theta})]$, and $\triangle_h [{\bm \Theta}(t), {\bm \Theta}'(t)] = \int \eta(\varphi_h) | \nabla_{\bm \Theta} \varphi_h \cdot  {\bm \Theta}'(t)|^2  \dd \x.$ For example, in Galerkin methods, one can let $\varphi_h(\x, t)$ be $\varphi_h(\x, t) = \sum_{i=1}^K \gamma_i(t) \psi_i(\x)$ with $\{ \psi_i \}_{i=1}^K$ being the set of basis functions. Then ${\bm \Theta}(t)$ is a vector given by ${\bm \Theta}(t) = (\gamma_1(t), \gamma_2(t), \ldots, \gamma_K(t))^{\rm T} \in \mathbb{R}^K$.

By employing the energetic variational approach in the semi-discrete level (\ref{discrte_En_1}), one can obtain an ODE system of ${\bm \Theta}(t)$. 
Particularly, in the linear response regime,  $\mathcal{D}_h ({\bm \Theta}(t), {\bm \Theta}'(t)) = \frac{1}{2} \int \eta(\varphi_h) |\nabla_{\bm \Theta} \varphi_h \cdot  {\bm \Theta}'(t)|^2 \dd \x$  is a quadratic function of ${\bm \Theta}'$.  The ODE system of ${\bm \Theta}(t)$ can then be written as
\begin{equation}
{\sf D} \left( {\bm \Theta} \right) {\bm \Theta}'(t)  =  -  \frac{\delta \mathcal{F}_h}{\delta {\bm \Theta}}~, % \quad   % \nabla_{{\bm \Theta}} \mathcal{F}_h(\bm{\Xi}(t))~,
\label{Eq:disc_ODE}
\end{equation}
where 
\begin{equation}
  \frac{\delta \mathcal{F}_h}{\delta {\bm \Theta}} = \frac{\delta \mathcal{F}}{\delta \varphi} \nabla_{\bm \Theta} \varphi_h~, \quad   {\sf D} \left( {\bm \Theta} \right) = \int \eta(\varphi_h)   (\nabla_{\bm \Theta} \varphi_h \otimes  \nabla_{\bm \Theta} \varphi_h )  \dd \x~.
\end{equation}
The ODE (\ref{Eq:disc_ODE}) is the same as the ODE (\ref{ODE_NGM}) in the neural Galerkin method \cite{bruna2022neural} for $L^2$ gradient flows, although the derivation is different.

Since (\ref{Eq:disc_ODE}) is a finite-dimensional gradient flow, one can then construct a minimizing movement scheme for ${\bm \Theta}$ \cite{Giorgi1992}: % or other variational integrators to (\ref{Eq:disc_ODE}). The minimizing movement 
finding ${\bm \Theta}^{n+1}$ such that
\begin{equation}\label{Min_Problem}
  \begin{aligned}
  & {\bm \Theta}^{n+1}  = \mathop{\arg\min}_{ {\bm \Theta} \in \mathcal{S}_{ad}^h} J_n( {\bm \Theta} )~, \quad J_n({\bm \Theta}) = \frac{ ( {\sf D}^n_*({\bm \Theta}  - {\bm \Theta}^{n}) ) \cdot ({\bm \Theta} - {\bm \Theta}^{n}) }{2 \tau} + \mathcal{F}_h ({\bm \Theta})~. \\
  % & \mathcal{S}_{ad}^h = \{ {\bm \Theta} \in \mathbb{R}^K ~|~ \det F_e ({\bm \Theta}) > 0, \quad e = 1, \ldots M.  \} \\
  \end{aligned}
\end{equation}
Here $\mathcal{S}_{ad}^h$ is the admissible set of ${\bm \Theta}$ inherited from the admissible set of $\varphi$, denoted by $\mathcal{S}$, and ${\sf D}^n_*$ is a constant matrix. A typical choice of ${\sf D}^n_* = {\sf D} ({\bm \Theta}^n)$. An advantage of this scheme is that
\begin{equation}
    \mathcal{F}_h ({\bm \Theta}^{n+1}) \leq   J_h^n({\bm \Theta}^{n+1}) \leq  J_h^n({\bm \Theta}^{n})   = \mathcal{F}_h ({\bm \Theta}^{n})~,
\end{equation}
if ${\sf D}^n_*$ is positive definite,  which guarantees the energy stability for the discrete free energy $\mathcal{F}_h ({\bm \Theta})$.  Moreover, by choosing a proper optimization method, we can assure that $\varphi^{n+1}_h$ stays in the admissible set $\mathcal{S}$.

Although neural networks can be used to construct $\varphi_h(\x; {\bm \Theta}(t))$, it might be expensive to compute $\nabla_{\bm \Theta} \varphi_h \otimes  \nabla_{\bm \Theta} \varphi_h$ in ${\sf D}({\bm \Theta})$. Moreover, ${\sf D}({\bm \Theta})$ is not a sparse matrix and requires a lot of computer memory to store when a deep neural network is used. 

To overcome these difficulties, we propose an alternative approach by introducing temporal discretization before spatial discretization. Let $\tau$ be the time step size. 
For the $L^2$-gradient flow (\ref{GD_1}), given $\varphi^n$, which represents the numerical solution at $t^n = n \tau$, one can obtain $\varphi^{n+1}$ by solving the following optimization problem: finding $\varphi^{n+1}$ in some admissible set $\mathcal{S}$ such that
\begin{equation}\label{Min_movement}
    \begin{aligned}
  & \varphi^{n+1} = \mathop{\arg\min}_{\varphi \in \mathcal{S}}   J^n (\varphi), \quad J^n (\varphi) =  \frac{1}{2 \tau} \int \eta(\varphi^n) | \varphi - \varphi^n  |^2 \dd \x + \mathcal{F}[\varphi]~.
    \end{aligned}
\end{equation}
Let $\varphi_h(\x; {\bm \Theta})$ be a finite-dimensional approximation to $\varphi$ with ${\bm \Theta} \in \mathbb{R}^K$ being the  parameter of the spatial discretization (e.g., weights of linear combination in a Galerkin approximation) yet to be determined, then the minimizing movement scheme (\ref{Min_movement}) can be written in terms of ${\bf \Theta}$: finding ${\bm \Theta}^{n+1}$ such that
\begin{equation}\label{Min_move_theta}
 {\bm \Theta}^{n+1} = \arg\min_{{\bm \Theta} \in \mathcal{S}_h} J_h^n ({\bm \Theta}), \quad J_h^n({\bm \Theta}) =  \frac{1}{2 \tau} \int \eta^n | \varphi_h(\x; {\bm \Theta}) - \varphi_h(\x; {\bm \Theta}^{n})|^2 \dd \x + \mathcal{F}_h [{\bm \Theta}]~.
\end{equation}
Here ${\bm \Theta}^{n}$ is the value of ${\bm \Theta}$ at time $t^n$. $\mathcal{F}_h [{\bm \Theta}] = \mathcal{F}[\varphi_h(\x; {\bm \Theta})]$, and $\eta^n = \eta( \varphi_h(\x; {\bm \Theta}^{n}) )$.

\begin{remark}
\label{rem:temp_space_dis}
The connection between the minimizing movement scheme (\ref{Min_move_theta}) derived by a temporal-then-spatial approach, and the minimizing movement scheme (\ref{Min_Problem}) derived by a spatial-then-temporal approach, can be shown with a direct calculation.
Indeed, according to the first-order necessary condition for optimality, an optimal solution to the minimization problem (\ref{Min_move_theta}) ${\bm \Theta}^{n+1}$ satisfies
\begin{equation}
\begin{aligned}
    \frac{\delta J^n_h({\bm \Theta})}{\delta {\bm \Theta}} \Big|_{ {\bm \Theta}^{n+1}} & = \frac{1}{\tau} \int \eta^n (\varphi_h(\x; {\bm \Theta}^{n+1} ) - \varphi_h(\x; {\bm \Theta}^n) ) \nabla_{\bm \Theta} \varphi_h  \Big|_{ {\bm \Theta}^{n+1}} \dd \x + \frac{\delta \mathcal{F}_h}{\delta {\bm \Theta}}  \Big|_{ {\bm \Theta}^{n+1}}  \\
    & = \int \eta^n \nabla_{\bm \Theta} \varphi_h  \Big|_{ {\bm \Theta}^{*}} \otimes \nabla_{\bm \Theta} \varphi_h  \Big|_{ {\bm \Theta}^{n+1}} \dd \x \frac{{\bm \Theta^{n+1}} - {\bm \Theta}^n}{\tau} + \frac{\delta \mathcal{F}_h}{\delta {\bm \Theta}}  \Big|_{ {\bm \Theta}^{n+1}} = 0    \\
  %   & = {\bm \Theta^{n+1}}} = \frac{1}{\tau} {\sf D}({\bm \Theta}^*) ({\bm \Theta}^{n+1} - {\bm \Theta}^n)      + \frac{\delta \mathcal{F}}{\delta {\bm \Theta}} = 0, \\
\end{aligned}
\end{equation}
for some ${\bm \Theta}^*$, where the second  equality follows the mean value theorem. In the case of Galerkin methods, as $\nabla_{\bm \Theta} \varphi_h$ is independent on ${\bm \Theta}$,
% Hence 
${{\bm \Theta}}^{n+1}$ is a solution to an implicit Euler scheme for the ODE (\ref{Eq:disc_ODE}) in which $\eta(\varphi_h)$ is treated explicitly.
\end{remark}

By choosing a certain neural network to approximate $\varphi$, denoted by $\varphi_{\rm NN}(\x; {\bm \Theta})$ (${\bm \Theta}$ is used to denote all parameters in the neural network), we can summarize the EVNN scheme for solving a $L^2$-gradient flow in Alg.~\ref{alg:L2_flow}.

\begin{algorithm}[H]
%\SetAlgoLined
% \SetAlgorithmName{Algorithm}{problem}{List of problems}
%\SetKwInOut{Input}{Input}\SetKwInOut{Output}{Output} %\LinesNumbered
 For a given initial condition $\varphi_0(\x)$, compute ${\bm \Theta}^0$ by solving
   \begin{equation}\label{Pre_train}
    {\bm \Theta}^0 = \mathop{\arg\min}_{\bm \Theta} \int_{\Omega} | \varphi_{\rm NN}(\x; {\bm \Theta}) - \varphi_0(\x)  |^2 \dd \x ~; % + \lambda \| {\bm w} \|^2~,
   \end{equation}
At each step, update ${\bm \Theta}^{n+1}$ by solving the optimization problem
 \begin{equation}
{\bm \Theta}^{n+1} = \mathop{\arg\min}_{\bm \Theta} \left( \frac{1}{2 \tau} \int_{\Omega} \eta^n  | \varphi_{\rm NN}(\x; {\bm \Theta}) - \varphi_{\rm NN}(\x, {\bm \Theta}^{n})|^2 \dd \x + \mathcal{F} [\varphi_{\rm NN}(\x; {\bm \Theta}) ]  \right).
\label{eq:innerloop}
 \end{equation}
 We have $\varphi_{\rm NN}(\x; {\bm \Theta}^n)$ as a numerical solution at time $t^n = n \tau$.  % {\color{blue} The last term can be viewed as a regularization term, since the first term may not be a good term (convex???) with respect to $\omega$}   

\caption{Numerical Algorithm for solving the $L^2$-gradient flow}
\label{alg:L2_flow}
\end{algorithm}

It can be noticed that both Eq.~(\ref{Pre_train}) and Eq.~(\ref{eq:innerloop}) involve integration in the computational domain $\Omega$. This integration is often computed by using a grid-based numerical quadrature or Monte--Carlo/Quasi--Monte--Carlo algorithms  \cite{rotskoff2020active}.
It is worth mentioning that due to the non-convex nature of the optimization problem and the error in estimating the integration, it might be difficult to find an optimal ${\bm \Theta}^{n+1}$ at each step. But since we start the optimization procedure with ${\bm \Theta}^n$,  we'll always be able to get a ${\bm \Theta}^{n+1}$ that lowers the discrete free energy at least on the training set. 
In practice, the optimization problem can be solved by either deterministic optimization algorithms, such as L-BFGS and gradient descent with Barzilai--Borwein step-size, or stochastic gradient descent algorithms, such as AdaGrad and Adam.

\begin{remark}
It is straightforward to incorporate other variational high-order temporal discretizations to solve the $L^2$-gradient flows. For example, a second-order accurate BDF2 scheme  can be reformulated as an optimization problem  
 \begin{equation}
 \begin{aligned}
{\bm \Theta}^n  = \mathop{\arg \min}_{\bm \Theta}  &\left( \frac{\eta}{ \tau} \int_{\Omega} \varphi_{\rm NN}(\x; {\bm \Theta}) - \varphi_{\rm NN}(\x, {\bm \Theta}^{n-1})|^2 \dd \x \right. \\ 
  & \left. - \frac{\eta}{4 \tau} \int_{\Omega}  |  \varphi_{\rm NN}(\x; {\bm \Theta}) - \varphi_{\rm NN}(\x, {\bm \Theta}^{n-2})|^2 \dd \x + \mathcal{F} [\varphi_{\rm NN}(\x; {\bm \Theta} ] \right)~. \\
  % + \gamma \|{\bm w} - {\bm w}^{n-1} \|^2 \right) , \\
\end{aligned}
 \end{equation}
 A modified Crank--Nicolson time-marching scheme can be reformulated as
  \begin{equation}
 \begin{aligned}
{\bm \Theta}^n  = \mathop{\arg \min}_{\bm \Theta}  &\left( \frac{\eta}{ \tau} \int_{\Omega}  |  \varphi_{\rm NN}(\x; {\bm \Theta}) - \varphi_{\rm NN}(\x, {\bm \Theta}^{n-1}) |^2 \dd \x + \mathcal{F} [\varphi_{\rm NN}(\x; {\bm \Theta}) ] \right.  \\
 & \left. \quad  + \langle  \nabla_{\Theta}\mathcal{F} [\varphi_{\rm NN}(\x; {\bm \Theta}^{n-1}) ],  {\bm \Theta} - {\bm \Theta}^{n-1} \rangle \right)~. \\
  % + \gamma \|{\bm w} - {\bm w}^{n-1} \|^2 \right) , \\
\end{aligned}
 \end{equation}
Here we assume that $\eta$ is a constant for simplicity.
 \end{remark}

\subsection{Lagrangian EVNN scheme for generalized diffusions}
% \textcolor{red}{add description here}

In this subsection, we show how to formulate an EVNN scheme in the Lagrangian frame of reference for generalized diffusions. 

As discussed previously, a generalized diffusion can be viewed as an $L^2$-gradient flow of the flow map $\x(\X, t)$ in the space of diffeomorphisms. Hence, the EVNN  scheme for the generalized diffusion can be formulated in terms of a minimizing movement scheme of the flow map  given by 
\begin{equation}\label{Flow_map}
   % \Phi^{n+1} = \mathop{\arg\min}_{\Phi \in {\rm Diff}}  \frac{1}{2 \tau} \int | \Phi(\X) - \Phi^n(\X) |^2    \rho_0(\X) \dd \X  + \mathcal{F} \left[\rho_0 \tilde{\circ} \Phi^{-1} (\x)  \right]~,
   \Phi^{n+1} = \mathop{\arg\min}_{\Phi \in {\rm Diff}}  \frac{1}{2 \tau} \int | \Phi(\X) - \Phi^n(\X) |^2    \rho_0(\X) \dd \X  + \mathcal{F} \left[  \Phi_{\verb|#|} \rho_0  \right]~,
\end{equation}
where $\Phi^{n}(\X)$ is a numerical approximation of the flow map $\x(\X, t)$ at $t^{n} = n \tau$, $\rho_0(\X)$ is the initial density. $\mathcal{F}[\rho]$
is the free energy for the generalized diffusion defined in Eq.~(\ref{F_diffusion_general}), and
$$ ( \Phi_{\verb|#|}\rho_0) (\x)  := \frac{\rho_0(\Phi^{-1}(\x))} {\det {\sf F} (\Phi^{-1}(\x))}~, \quad {\rm Diff} = \{ \Phi: \mathbb{R}^d \rightarrow \mathbb{R}^d ~|~ \Phi~\text{is a diffeomorphism}  \}~.$$ 
% $\tilde{\circ}$ is a kinematics-preserving composition.
One can parameterize $\Phi: \mathbb{R}^d \rightarrow \mathbb{R}^d$ by a suitable neural network. The remaining procedure is almost the same as that of the last subsection. However, it is often difficult to solve this optimization problem directly, and one might need to build a large neural network to approximate $\Phi^{n+1}$ when $n$ is large.

To overcome this difficulty, we proposed an alternative approach, Instead of seeking an optimal map $\Phi^{n+1}$ between $\rho^0$ and $\rho^{n+1}$, we seek an optimal $\Psi^{n+1}$ between $\rho^n$ and $\rho^{n+1}$. More precisely.
we assume that 
$$\Phi^{n+1} = \Psi^{n+1} \circ \Psi^n \circ \Psi^{n-1} \ldots \circ \Psi^1~.$$  Given $\rho^n$, one can compute $\Psi^{n+1}$ by solving the following optimization problem
\begin{equation}\label{Scheme_LG1}
    \Psi^{n+1} = \mathop{\arg\min}_{\Psi \in {\rm Diff}}  \frac{1}{2 \tau} \int | \Psi(\x) - \x |^2    \rho^n(\x) \dd \x  + \mathcal{F} [ \Psi_{\verb|#|} \rho^n ]~.
\end{equation}
The corresponding $\rho^{n+1}$ can then be computed through
$ \rho^{n+1}(\x) = (\Psi^{n+1}_{\verb|#|} \rho^n ) (\x)$. An advantage of this approach is that we only need a small size of neural network to approximate $\Psi^{n+1}$ at each time step when $\tau$ is small.

\begin{remark}
The scheme (\ref{Scheme_LG1}) can be viewed as a Lagrangian realization of the JKO scheme (\ref{JKO}) for the Wasserstein gradient flow, although it is developed based on the $L^2$-gradient flow structure in the space of diffeomorphism. According to the Benamou--Brenier formulation {\rm \cite{benamou2000computational}}, the Wasserstein distance between two probability densities $\rho_1$ and $\rho_2$ can be computed by solving the optimization problem 
\begin{equation}\label{CFD_Wd}
W_2(\rho_1, \rho_2)^2 = \mathop{\min}_{(\rho, \uvec) \in \mathcal{S}} \int_{0}^1 \int \rho |\uvec|^2 \dd \x \dd t~,
\end{equation}
where the admissible set of $(\rho, \uvec)$ is given by
\begin{equation}
\mathcal{S} = \{ (\rho, \uvec) ~|~  \rho_t + \nabla \cdot (\rho \uvec) = 0~,\quad \rho(\x, 0) = \rho_1~, \quad \rho(\x, 1) = \rho_2 \}~.
\end{equation}
Hence, one can solve the JKO scheme by solving
\begin{equation}\label{JKO_OC}
  \begin{aligned}
    &  (\uvec^*, \hat{\rho}^*) = \mathop{\arg\min}_{(\uvec, \hat{\rho})} \frac{1}{2} \int_{0}^{\tau} \int_{\Omega} \hat{\rho} |\uvec|^2 \dd \x \dd t + \mathcal{F}[\hat{\rho}(\tau)] \\
    & \text{s.t.} \quad    \pp_t \hat{\rho} + \nabla \cdot (\hat{\rho} \uvec) = 0, \quad \hat{\rho}(0) = \rho^{n}. \\
    \end{aligned}
\end{equation}
and letting $\rho^{n+1} = \hat{\rho}^*(\tau)$. In Lagrangian coordinates, (\ref{JKO_OC}) is equivalent to
\begin{equation}\label{JKO_OC_Lag}
   \x^*(\X, t) = \mathop{\arg\min}_{\x(\X, t)} \frac{1}{2} \int_{0}^{\tau} \int_{\Omega} \rho_0(\X) |\x_t(\X, t)|^2 \dd \X \dd t + \mathcal{F} (\hat{\rho}(\x, \tau)), 
\end{equation}
where $\rho_0 = \rho^n$ and $\hat{\rho}(\x(\X, \tau), \tau) = \rho_0(\X) / \det {\sf F}(\X, \tau)$. By taking variation of (\ref{JKO_OC_Lag}) with respect to $\x(\X, t)$, one can show that the optimal condition is $\x_{tt}(\X, t) = 0$ for $t \in (0, \tau)$, which indicates that $\x(\X, t) = t (\Psi(\X) - X) / \tau + X$ if $\x(\X, \tau) = \Psi(\X)$.  
Hence, if $\Psi^*$ is the optimal solution of (\ref{Scheme_LG1}), then  $\x^*(\X, t) = t (\Psi^*(\X) - \X) / \tau + \X$ is the optimal solution of  (\ref{JKO_OC_Lag}).

\end{remark}

\begin{remark}
If $\eta(\rho) \neq \rho$, we can formulate the optimization problem (\ref{Scheme_LG1}) as
\begin{equation}
    \Psi^{n+1} = \mathop{\arg\min}_{\Psi \in {\rm Diff}}  \frac{1}{2 \tau} \int | \Psi(\x) - \x |^2    \eta( \rho^n(\x)) \dd \x  + \mathcal{F} [ \Psi_{\verb|#|} \rho^n  ]~.
\end{equation}
by treating $\eta$ explicit. A subtle fact is that $\det {\sf F}^n = 1$ since we always start with an identity map.
\end{remark}
 
The numerical algorithm for solving the generalized diffusion is summarized in Alg.~\ref{alg:gen_diff}.
\begin{algorithm}[H]
\begin{itemize}
\item Given $\{ \x_i^n \}_{i=1}^N$ and the densities $\rho^n_i$ at $t = n \tau$ and $\x_i^n$.

\item Find $\Psi^{n+1}(\x): \mathbb{R}^d \rightarrow \mathbb{R}^d$, by solving the optimization problem (\ref{Scheme_LG1}).
 To guarantee energy stability, we should take $\Psi$ as an approximation to an identity map initially when solving the optimization problem (\ref{Scheme_LG1}).

\item After obtaining $\Psi^{n+1}$, updates $\{ \x_i^{n+1} \}_{i=1}^N$ and $\rho^{n+1}$ by
\begin{equation}
\begin{aligned}
&   \x_i^{n+1} = \Psi^{n+1} (\x_i^{n})~, \quad \rho^{n+1}_i = \frac {\rho^n_i }{\det ( \nabla \Psi^{n+1} (\x_i^n) )}~. \\
\end{aligned}
\end{equation}
% \item {\bf Resampling} We might need resample. Generate a new set of samples $\{  \tilde{\x}_i^{k+1} \}_{i=1}^N$, compute $\rho^{k+1}$ in this new set of sample. And using $\{ \tilde{\x}_i^{k+1}  \}_{i=1}^{N}$ in (\ref{Up_PHI}) to update $\Phi^{k+1}$.
\end{itemize}
\caption{Numerical Algorithm for solving the generalized diffusion}
\label{alg:gen_diff}
\end{algorithm}

 The next question is how to accurately evaluate the numerical integration in (\ref{Scheme_LG1}). Let $\x_i^n = \Phi^n (\x_i^0)$, for the general free energy (\ref{F_diffusion_general}),  %the loss function 
 one way to evaluate the integrations in Eq.~(\ref{Scheme_LG1}) is using
\begin{equation}
\begin{aligned}
  & J(\Psi) =  \textstyle \frac{1}{2 \tau}  \sum_{i=1}^N \rho^n_i \| \Psi(\x_i^n) - \x_i^n  \|^2  |\Omega_i^n|  
%  & +  \sum_{i=1}^N \left( \rho^n_i ( \ln (\rho^n_i) -  \ln (\det ( \nabla \Psi (\x_i^k) )) ) + \rho_i^n V(\Psi(\x_i^n)) \right)  |\Omega_i^n| ~, \\
%  & +  \sum_{i=1}^N \left( \omega(\rho_i^n / \det ( \nabla \Psi (\x_i^k) ) ) \det (\nabla \Psi (\x_i^n) ) + V(\Psi(\x_i^n)) \rho_i^n   \right)   |\Omega_i^n| ~, \\
%  & +  \frac{|\Omega|}{N}  \sum_{i=1}^N \left( \rho^k(\x_i^k) \ln (\rho^k (\x_i^k)) - \rho^k(\x_i^k) \ln (\det ( \nabla \Phi (\x_i^k)  ))  + \rho^k(\x_i^k) V(\Phi(\x_i^k))   \right) \\
  +  \sum_{i=1}^N \left( \omega\left( \frac{\rho_i^n} { \det ( \nabla \Psi (\x_i^n) )} \right) \det (\nabla \Psi (\x_i^n) )  \right)   |\Omega_i^n|\\
  & \quad \quad  + \textstyle \sum_{i=1}^N V(\Psi(\x_i^n)) \rho_i^n |\Omega_i^n| +   \frac{1}{2} \sum_{i, j=1}^N K( \Psi(\x_i), \Psi(\x_j)) \rho_i^n  \rho_j^n |\Omega_j^n| |\Omega_i^n|  ~, \\
 \end{aligned}
 \end{equation}
 where $|\Omega_i^n|$ is the volume of the Voronoi cells associated with the set of points $\{ \x_i^n \}$, $\rho^n_i$ stands for $\rho^n(\x_i^n)$. 
Here the numerical integration is computed through a piece-wisely constant reconstruction of $\rho^n$ based on its values at $\{ \Phi^n(\x_i^0) \}_{i=1}^N$. % {\color{blue}So $\rho_i^n$ should be understood as volume average.} 
However, it is not straightforward to compute $|\Omega_i^n|$, particularly for high dimensional cases.
In the current study, we assume that the initial samples are drawn from $\rho_0$, one can roughly assume that $\{ \x_i^n \}$ follows the distribution $\rho^n$, then according to the Monte--Carlo approach, the numerical integration can be evaluated as
\begin{equation}
\begin{aligned}
  J(\Psi) & = \textstyle  \frac{1}{2 \tau} \left( \frac{1}{N} \sum_{i=1}^N   \| \Psi(\x_i^n) - \x_i^n  \|^2 \right)  \\
%  & +  \frac{1}{N} \sum_{i=1}^N \bigg( \big( \ln (\rho^n_i) -  \ln (\det ( \nabla \Psi (\x_i^n)  )) +  V(\Psi(\x_i^n)) \big) \bigg) \\
%  & +  \frac{|\Omega|}{N}  \sum_{i=1}^N \left( \rho^k(\x_i^k) \ln (\rho^k (\x_i^k)) - \rho^k(\x_i^k) \ln (\det ( \nabla \Phi (\x_i^k)  ))  + \rho^k(\x_i^k) V(\Phi(\x_i^k))   \right) \\
  & \textstyle  + \frac{1}{N} \sum_{i=1}^N \left( f_{\omega} \left(  \frac{\rho_i^n }{ \det ( \nabla \Psi (\x_i^n) )} \right)  + V(\Psi(\x_i^n))   \right)   +  \frac{1}{2 N^2} \sum_{i, j=1}^N K \left( \Psi(\x_i), \Psi(\x_j)\right)  ~, \\
 \end{aligned}
 \end{equation}
 where $f_{\omega} (\rho) = \omega(\rho) / \rho.$ The proposed numerical method can be further improved if one can evaluate the integration more accurately, i.e., have an efficient way to estimate $|\Omega_i^n|$.  Alternatively, 
 as an advantage of the neural network-based algorithm, we can get $\rho^n(x) = \rho_0 \tilde{\circ} (\Phi^n)^{-1} (\x), ~ \forall \x$, which enables us to estimate $J(\Psi)$ by resampling. More precisely, assume $\mu$ is a distribution that is easy to sample
 \begin{equation}
 \begin{aligned}
     & J (\Psi) = 
     \textstyle   \frac{1}{2 \tau} \left( \frac{1}{N} \sum_{i=1}^N   \| \Psi(\x_i^n) - \x_i^n  \|^2 \right) \frac{\rho_i^n}{\mu_i^n}  + \frac{1}{N} \sum_{i=1}^N \left( f_{\omega} \left(  \frac{\rho_i^n }{ \det ( \nabla \Psi (\x_i^n) )} \right)  \right) \frac{\rho_i^n}{\mu_i^n}  \\
 %  & +  \frac{1}{N} \sum_{i=1}^N \bigg( \big( \ln (\rho^n_i) -  \ln (\det ( \nabla \Psi (\x_i^n)  )) +  V(\Psi(\x_i^n)) \big) \bigg) \\
 %  & +  \frac{|\Omega|}{N}  \sum_{i=1}^N \left( \rho^k(\x_i^k) \ln (\rho^k (\x_i^k)) - \rho^k(\x_i^k) \ln (\det ( \nabla \Phi (\x_i^k)  ))  + \rho^k(\x_i^k) V(\Phi(\x_i^k))   \right) \\
   & \quad \textstyle +  V(\Psi(\x_i^n))  \frac{\rho_i^n}{\mu_i^n}   +  \frac{1}{2 N^2} \sum_{i, j=1}^N K \left( \Psi(\x_i), \Psi(\x_j)\right) \frac{\rho_i^n \rho_j^n}{\mu_i^n \mu_j^n}  ~, \\
  \end{aligned}
 \end{equation}
 where $\x_i \sim \mu$, $\rho_i^n = \rho^n (x_i)$. We will explore this resampling approach in future work.

The remaining question is how to parameterize a diffeomorphism using neural networks. This is discussed in detail in the next subsection.

\subsection{Neural network architectures}

In principle, the proposed numerical framework is independent of the choice of neural network architectures. However, different neural network architectures may lead to different numerical performances, arising from a balance of approximation (representation power), optimization, and generalization. In this subsection, we briefly discuss several neural network architectures that we use in the numerical experiments.

\subsubsection{Neural network architectures for Eulerian methods} 

For Eulerian methods, one can construct a neural network to approximate the unknown function $f: \mathbb{R}^d \rightarrow \mathbb{R}$. Shallow neural networks (two-layer neural networks) approximate $f$ by functions of the form
\begin{equation}\label{SN_ap}
\begin{aligned}
  f(\x; {\bm \Theta}) & = \sum_{i = 1}^N \alpha_i \sigma(\bm{\omega}_i \cdot \x + b_i) + \alpha_0 = \bm{\alpha} \cdot \sigma ({\bm W} \x + \bm{b}) + \alpha_0 ~, \\
 \end{aligned}
\end{equation}
where $\sigma(.): \mathbb{R} \rightarrow \mathbb{R}$ is a fixed nonlinear activation function, $N$ is the number of hidden nodes (neurons), and ${\bm \Theta} = \{ \alpha_i,{\bm \omega}_i, b_i \}$ are the NN parameters to be identified. Typical choices of activation functions include the ReLU $\sigma(x) = \max(x, 0)$, the sigmoid $\sigma(x) = 1 / (1 + e^{-2x})$ and the hyperbolic tangent function $\tanh(x)$.% ${\bm \alpha} = (\alpha_1, \alpha_2, \ldots, \alpha_N)^{\rm T}$
A DNN can be viewed as a network composed of many hidden layers. More precisely, a DNN with $L$ hidden layers represents a function $f: \mathbb{R}^{d} \rightarrow \mathbb{R}$ by \cite{shen2019nonlinear} 
\begin{equation}
f(\x; {\bm \Theta}) = g \circ T^{L} \circ T^{L-1} \circ \ldots \circ T^{(1)} (\x)~, 
\end{equation}
where $g(\z) = \sum_{i=1}^{N_L} \gamma_i z_i + \gamma_0$ is a linear map from $\mathbb{R}^{N_l}$ to $\mathbb{R}$, $T^{(l)}$ is a nonlinear map from $\mathbb{R}^{N_{l-1}}$ to $\mathbb{R}^{N_{l}}$, and $N_0 = d$. The nonlinear map $T^{l}$ takes the form
\begin{equation}\label{def_Tl}
T^{(l)} (\x_{l-1})  = \sigma ({\bm W}_{l} \x_{l-1} + \bm{b}_l)~,
\end{equation}
where ${\bm W}_l \in \mathbb{R}^{N_{l-1} \times N_l}$, $\bm{b}_l \in \mathbb{R}^l$ and $\sigma(\cdot)$ is a nonlinear activation function that acts component-wisely on vector-valued inputs.
%{\color{red} Why not use shallow network }
% 

Another widely used class of DNN model is residual neural network (ResNet).  A typical $K$-block ResNet approximates an unknown function $f(\x): \mathbb{R}^d \rightarrow \mathbb{R}$ by
\begin{equation}
    f_K(\x; {\bm \Theta}) = g(\z_K(\x))~,
\end{equation}
where $g(\z) = \sum_{i=1}^N \gamma_i z_i + \gamma_0$ is a linear map from $\mathbb{R}^N \rightarrow \mathbb{R}$ and $z_K(\x): \mathbb{R}^d \rightarrow \mathbb{R}^N$ is a nonlinear map defined through
\begin{equation}\label{def_zL}
    \begin{aligned}
     \z_0 =  {\bm V} \x~, \quad
        %\z_i = \sigma_2 ( {\bm \alpha}_i  \sigma_1 ( {\bm W}_i \z_{i-1} + {\bm b}_i ) + \z_{i-1} )~, 
         \z_k = \sigma_2 ( {\bm \alpha}_k  T^{L_k}_k \circ T^{L-1}_k \circ \ldots \circ T^{(1)}_k (\z_{k-1})  + \z_{k-1} )
        \quad k = 1, 2, \ldots, K~.
    \end{aligned}
\end{equation}
Here $L_i$ is the number of fully connected layer in $i$-the block, $T^{(l)}_k$ is the same nonlinear map defined in (\ref{def_Tl}) with ${\bm W}_l^k \in \mathbb{R}^{M \times N}$, ${\bm b}_l^k \in \mathbb{R}^M$, ${\bm V} \in \mathbb{R}^{N \times d}$, ${\bm \alpha}_i \in \mathbb{R}^{N \times M}$ and $\sigma_i(\cdot)$ is an element-wise activation function. The model parameters are ${\bm \Theta} = \{ {\bm \gamma}, {\bm \alpha}_i, {\bm W}_l^k, {\bm b}_l^k, {\bm V}   \}$. 
The original ResNet \cite{he2016deep} takes $\sigma_2$ as a nonlinear activation function such as ReLU. 
Later studies indicate that one can also take $\sigma_2$ as the identity function \cite{weinan2021algorithms, he2016identity}. Then at an infinite length, i.e., $L \rightarrow \infty$, (\ref{def_zL}) corresponds to the ODE 
\begin{equation}
    \frac{\dd \z}{\dd t} = f(\z)~, \quad \z_0 = \x~.
\end{equation}
Compared with fully connected deep neural network models which may suffer from numerical instabilities in the form of exploding or vanishing gradients \cite{EMaWu2020, Hanin2018WhichNN, Kolen_Kremer01},  very deep ResNet can be constructed to avoid these issues.

Given that the proposed numerical scheme employs neural networks with time-dependent parameters to approximate the solution of a gradient flow, there is no need to employ a deep neural network. In all numerical experiments for $L^2$-gradient flows, we utilize ResNet with only a few blocks. The detailed neural network settings for each numerical experiment will be described in the next section. We'll compare the numerical performance of different neural network architectures in future work.

\subsubsection{Neural network architectures for Lagrangian methods}

The proposed Lagrangian method seeks for a neural network to approximate a diffeomorphism from $\mathbb{R}^d$ to $\mathbb{R}^d$.  This task involves two main challenges: one is ensuring that the map is a diffeomorphism, and the other is efficiently and robustly computing the deformation tensor ${\sf F}$ and its determinant $\det {\sf F}$.
Fortunately, various neural network architectures have been proposed to approximate a transport map.  Examples include planar flows \cite{rezende2015variational}, auto-regressive flows \cite{kingma2016improved}, continuous-time diffusive flow \cite{tabak2010density}, neural spline flow \cite{durkan2019neural}, and
convex potential flow \cite{huang2020convex}.

One way to construct a neural network $\mathbb{R}^d \rightarrow \mathbb{R}^d$ for approximating a flow map is to use a planar flow. A $K$-layer planar flow is defined by  $T = T_K \circ \ldots   T_1 \circ  T_0$, where $T_k: \mathbb{R}^d \rightarrow \mathbb{R}^d$ is given by
\begin{equation}
    \x^{k+1} = T_k (\x^{k}) = \x^{k} +  {\bm u}_k h ({\bm w}^{\rm T}_k \x^{k} + b_k )~.
\end{equation}
Here, ${\bm w}_k, {\bm u}_k \in \mathbb{R}^d$, $b_k \in \mathbb{R}$, and $h$ is a smooth, element-wise, non-linear activation function such as $\tanh$. Direct computation shows that
\begin{equation}
J_k = \det (\nabla T_k) = 1 + h'({\bm w}^{\rm T}_k \x^{k} + b_k ) {\bm u}^{\rm T}_k {\bm w}_k ~.
\end{equation}
Clearly, $T_k$ is a diffeomorphism if $ h'({\bm w}^{\rm T}_k \x^{k} + b_k ) {\bm u}^{\rm T}_k  {\bm w}_k  < 1$ for $\forall \x_k$.
The determinant of the transport map can be computed as $\det (\nabla T) = J_K J_{K-1} \ldots J_0$. and we have $\varphi(T(\x)) = \frac{\varphi^n(\x) }{ \det (\nabla T)}$.

One limitation of planar flow is its potential lack of expressive power. Another commonly employed neural network architecture for approximating a flow map is the
convex potential flow \cite{huang2020convex}, which defines a diffeomorphism via the gradient of a strictly convex function that is parameterized by an Input Convex Neural Network (ICNN) \cite{amos2017input}.
% parameterizes an unknown convex function which guarantees the existence of diffeomorphism.
A fully connected, $K-$layer ICNN can be written as
\begin{equation}
\z_{l+1} = \sigma_l ( {\bm W}_l \z_l + \bm{A}_l \X + \bm{b}_l)~, \quad l = 0, 1, \ldots K-1~, %\quad f(\X; \theta) = z_l,
\end{equation}
where $\z_0 = 0$, $\bm{W}_0 = {\sf 0}$, ${\bm \Theta} = \{ \bm{W}_l, \bm{A}_l, \bm{b}_l  \}$ are the parameters to be determined, and $\sigma_l$ is the nonlinear activation function. As proved in \cite{amos2017input}, if all entries $\bm{W}_l$ are non-negative, and all $\sigma_l$ are convex and non-decreasing, then $f(\X; {\bm \Theta})$ is a convex function with respect to $\X$. Hence $\nabla_{\X} f(\X; {\bm \Theta})$ provides a parametric approximation to a flow map. In the current study, we adopt the concept of the convex potential flow to develop the Lagrangian EVNN method. We'll explore other types of neural network architectures in future work.

\begin{remark}
It is worth mentioning that the gradient of a convex function $f$ only defines a subspace of diffeomorphism, and it is unclear whether the optimal solution of (\ref{Scheme_LG1}) belongs to this subspace.
\end{remark}

\section{Numerical Experiments}
\label{sec:simu}
In this section, we test the proposed EVNN methods for various $L^2$-gradient flows and generalized diffusions. To evaluate the accuracy of different methods, we define the $l^2$-error $$\left( \textstyle \frac{1}{N} \sum_{i=1}^N |\varphi_{\rm NN}(\x_i) - \varphi_{\rm ref}(\x_i)|^2 \right)^{1/2}$$
%$\int_{\Omega} |\varphi_{\rm ref} - \varphi_{\rm NN}|^2 \dd \x$,
and the relative $l_2$-error $$\left( \textstyle \sum_{i=1}^N |\varphi_{\rm NN}(\x_i) - \varphi_{\rm ref}(\x_i) |^2 /  (\textstyle \sum_{i=1}^N |\varphi_{\rm ref}(\x_i))|^2 \right)^{1/2}$$ between the NN solution $\varphi_{\rm NN}$ and the corresponding reference solution $\varphi_{\rm ref}$ on a set of test samples $\{ \x_i \}_{i=1}^N$. The reference solution $\varphi_{\rm ref}(\x)$ is either an analytic solution or a numerical solution obtained by a traditional numerical method. In some numerical examples, we also plot point-wise absolute errors $|\varphi_{\rm ref} - \varphi_{\rm NN}|$ at the grid points.

\subsection{Poisson equations}
Although the proposed method is designed for evolutionary equations, it can also be used to compute equilibrium solutions for elliptic problems. In this subsection, we compare the performance of the EVNN method with two classical neural network-based algorithms, PINN and DRM, in the context of solving Poisson equations.

We first consider a 2D Poisson equation with a Dirichlet boundary condition
\begin{equation}
\label{eq:poisson}
    \begin{cases}
    & -\Delta u(\x) = f(\x)~, \quad  \x \in \Omega \subset \mathbb{R}^d~, \\
  & u(\x) = g(x)~,\quad \quad \quad \quad \x \in \partial\Omega~. \\
    \end{cases}
\end{equation}
Since EVNN is developed for evolution equations, we solve the following $L^2-$gradient flow
\begin{equation}
\frac{\dd}{\dd t} \left( \int_{\Omega} \frac{1}{2} |\nabla u|^2  - f(\x) u(\x) \dd \x + \lambda \int_{\pp \Omega} |u(\x) - g(\x)|^2  \dd S   \right)= - \int_{\Omega} |u_t|^2 \dd \x~,
\end{equation}
to get a solution of the Poisson equation (\ref{eq:poisson}).
Here, $\lambda \int_{\pp \Omega} |u - g(\x
)|^2  \dd S$ is the surface energy that enforces the Dirichlet boundary condition. The corresponding gradient flow equation is given by
\begin{equation}
  u_t = \Delta u (x)  + f(\x),
\end{equation}
subject to a Robin boundary condition $u(\x) - g(\x) = \frac{1}{ 2 \lambda} \frac{\pp u}{\pp {\bm n}}$, where ${\bm n}$ is the outer normal of $\pp \Omega$. One can recover the original Dirichlet boundary condition by letting $\lambda \rightarrow \infty$.
Such a penalty approach is also used in PINN and DRM, and we take $\lambda = 500$ in all numerical experiments below.

We consider the following two cases:
\begin{itemize}
    \item Case 1: $\Omega = (0, \pi) \times (-\pi/2, \pi/2)$, $f(\x) = 2 \sin x \cos y$ and $g(\x) = 0$. The exact solution is $u(\x) = \sin x \cos y$.
    \item Case 2: $\Omega = \{ (x, y) ~|~ | \x | \le 1 \}$, $f(\x) = \frac{\pi^2}{4} \sin ( \frac{\pi}{2} (1 - |\x|)) + \frac{\pi}{2 |\x|} \cos (\frac{\pi}{2} (1 - |\x|)) $ and $g(\x) = 0$. The exact solution is $u(\x) = \sin (\frac{\pi}{2} (1 - |\x|))$.
\end{itemize}
We employ a 1-block ResNet with 20 hidden nodes and one hidden layer for all cases. The total number of parameters is 501. We apply Xavier Initialization \cite{glorot2010understanding} to initialize the weights of neural networks in all cases. To evaluate the integration in all methods, we generate $2500$ samples in $\Omega$ using a Latin hypercube sampling (LHS) and $50$ samples on each boundary of $\pp \Omega$ using a 1D uniform sampling for case 1. For case 2, we generate $2500$ samples in $(-1, 1)^2$ using LHS, but only use the samples satisfies $x^2 + y^2 < 1$ as training samples. Additionally, we generate 200 training samples $(\cos \theta_i, \sin \theta_i)_{i=1}^{200}$ on the boundary, with  $\{ \theta_i \}_{i=1}^{200}$ being generated by a uniform distribution on $(0, 2 \pi)$.
For PINN and DRM, we use Adam to minimize the loss function and use a different set of samples at each iteration. For EVNN, we use a different set of samples at each time step and employ an L-BFGS to solve the optimization problem (\ref{eq:innerloop}).

\begin{figure}[!htb]
  \includegraphics[width = \linewidth]{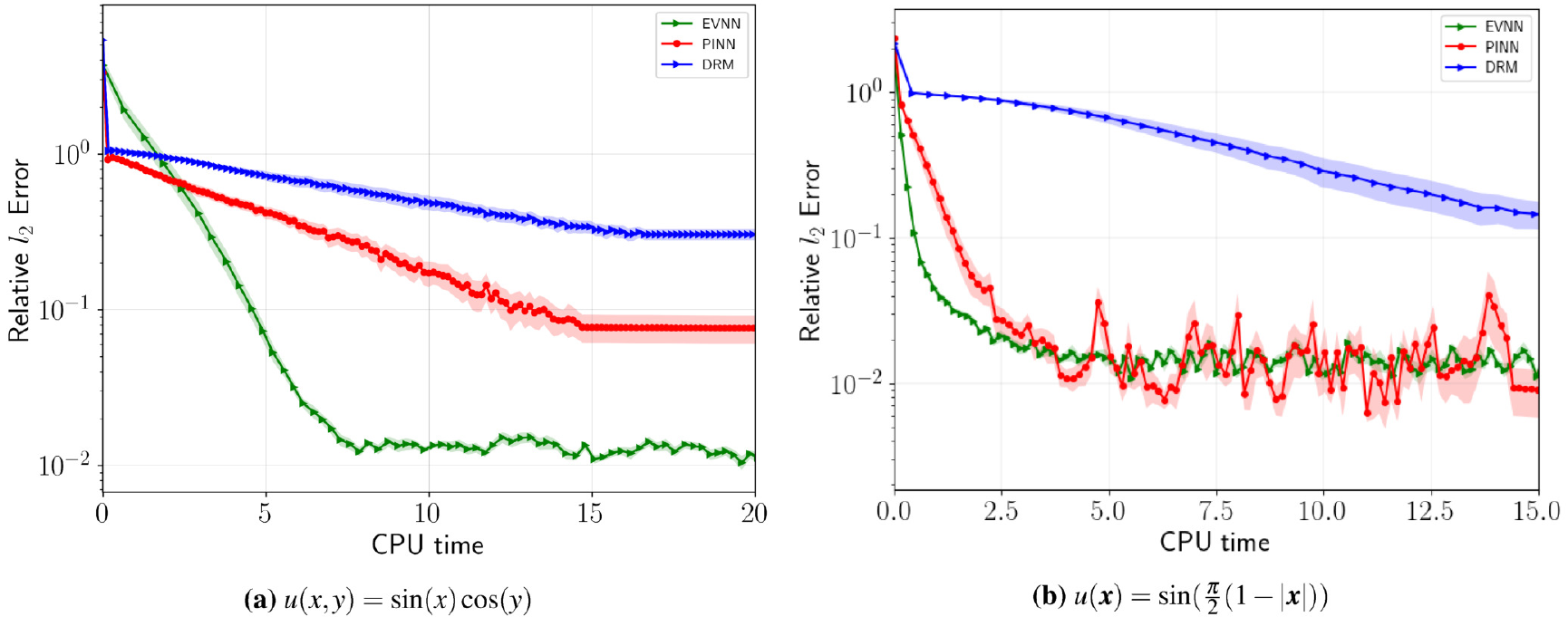}
  \caption{ Relative $l^2$ error with respect to CPU time for different methods for 2D Poisson equations. 
  }\label{fig:poisson2d}
  \end{figure}

Due to randomness arising from the stochastic initialization of the neural network and the sampling process, we repeated each method 10 times and plotted the mean and standard error of the relative $l_2$-errors.
Fig. \ref{fig:poisson2d}(a) shows the relative $l_2$-errors of different methods with respect to the CPU time for both cases. The training was executed on a MacBook Pro equipped with an M2 chip. In Case 1, the relative $l_2$ error of each method is evaluated on a test set comprising a uniform grid of $101 \times 101$ points within the domain $(0, 1)^2$. For Case 2, the test set is defined as the collection of points $\{(x_i, y_j) | x_i^2 + y_j^2 < 1 \}$, where $(x_i, y_j)$ forms a uniform $201 \times 201$ grid on $(-1, 1)^2$. The CPU time is computed by averaging the time of 10 trials at $n$-th iteration. The numerical results clearly show that the proposed method is more efficient than both PINN and DRM. It significantly enhances the efficiency and test accuracy of DRM through the introduction of time relaxation. While PINN may achieve better test accuracy in Case 2, it exhibits a larger standard error in both cases.

Next, we consider a high-dimensional Poisson equation
\begin{equation}
% \begin{cases}
    - \Delta u = f(\x)~, \quad \x \in \Omega =  (-1, 1)^d~, \\
%    & \frac{\pp u}{\pp n} |_{\pp \Omega} = 0 ~.\\
% \end{cases}
\label{eq:Hi_Poisson}
\end{equation}
with a homogeneous Neumann boundary condition $\frac{\pp u}{\pp n} |_{\pp \Omega} = 0$. We take $f(\x) = \pi^2 \sum_{k=1}^d \cos (\pi x_k)$.
Similar numerical examples are tested in \cite{weinan2018deep, lu2021priori}. The exact solution to this high-dimensional Poisson equation is $u(\x) = \sum_{k=1}^d \cos(\pi x_k)$. 
% The energy-dissipation law for
Following \cite{lu2021priori}, solve an $L^2$-gradient flow associated with the free energy
\begin{equation}
    \mathcal{F}[u] = \textstyle \int_{(-1, 1)^d} \frac{1}{2} |\nabla u|^2 - fu ~\dd \x + \lambda \left(\int_{(-1, 1)^d} u ~\dd \x \right)^2,
\end{equation} 
where the last term enforces $\int_{[-1, 1]^d} u ~\dd \x = 0$. We take $\lambda = 0.5$ as in \cite{lu2021priori} in all numerical tests.

\begin{table}[ht]
  \caption{Relative $l^2$-error of high-dimensional Poisson equation with different settings in different dimensions. The first column shows the number of residual blocks, the number of nodes in each fully connected layer, and the number of samples.}
  \centering 
  \begin{tabular}{c | c | c | c | c} 
  \hline
   Setting  & $d = 4$ & $d = 8$ & $d = 16$ & $d = 32$\\ [0.5ex]
  \hline
   (2, 10, 1000) & 0.024  & 0.046  & 0.623  & 0.867 \\
  \hline
   (3, 60, 1000) & 0.042  & 0.048 & 0.077  & 0.117  \\
  \hline
   (3, 60, 10000) & 0.018  & 0.022  & 0.036 & 0.077 \\
  \hline
  
  \end{tabular}
  \label{table:highpoi} 
  \end{table}

Since it is difficult to get an accurate estimation of high-dimensional integration, we choose Adam to minimize (\ref{eq:innerloop}) at each time step and use a different set of samples at each Adam step. This approach allows us to explore the high dimensional space with a low computational cost. The samples are drawn using LHS.
Table \ref{table:highpoi} shows the final relative $l^2$-error with different neural network settings in different spatial dimensions. The relative $l^2$-error is evaluated on a test set that comprises 40000 samples generated by LHS.
The model is trained with $200$ outer iterations ($\tau = 0.01$) and at most $200$ iterations for inner optimization (\ref{eq:innerloop}). An early stop criterion was applied if the $l_2$-norm of the model parameters between two consecutive epochs is less than $10^{-6}$. The first column of Table~\ref{table:highpoi} specifies the numerical setting as follows: the entry $(2, 10)$ means a neural network with 2 residual blocks, each with two fully connected layers having 10 nodes is used, and $1000$ represents the number of samples drawn in each epoch. As we can see, the proposed method can achieve comparable results with results reported in similar examples in previous work by DRM \cite{weinan2018deep, lu2021priori}. It can be observed that increasing the width of the neural network improves test accuracy in high dimensions, as it enhances the network's expressive power. Furthermore, increasing the number of training samples also significantly improves test accuracy in high dimensions. We'll explore the effects of $\tau$, the number of samples, and neural network architecture for the high-dimensional Poisson equation in future work.

\subsection{$L^2$-gradient flow}
In this subsection, we apply the EVNN scheme to two $L^2$-gradient flows to demonstrate its energy stability and numerical accuracy.
\subsubsection{Heat equation}
We first consider an initial-boundary value problem of a heat equation:
\begin{equation}
    \begin{cases}
        & u_t = \Delta u({\bf x}), \quad {\bf x} \in \Omega =(0, 2)^2~,\quad t \in (0, T]~,\\
        & u(\x, t) = 0, \quad \x \in \partial\Omega~, \quad t \in (0, T]~,\\
        & u(\x, 0) = \sin \left( \frac{\pi}{2} x_1 \right) \sin \left( \frac{\pi}{2} x_2 \right)~, \quad \x \in \Omega = (0, 2)^2~,\\
    \end{cases}
\end{equation}
which can be interpreted as an $L^2-$gradient flow satisfying the energy-dissipation law 
\begin{equation}
   \textstyle \frac{\dd}{\dd t} \left( \int_{\Omega}\frac{1}{2}\left|\nabla u\right|^2 \dd \x + \lambda \int_{\pp \Omega} |u|^2 \dd S \right) = -\int \left|u_t\right|^2 \dd \x~.
\end{equation}
Again, a surface energy term is added to enforce the Dirichlet boundary condition.

\begin{figure}[!tbh]
\includegraphics[width = 0.95 \linewidth]{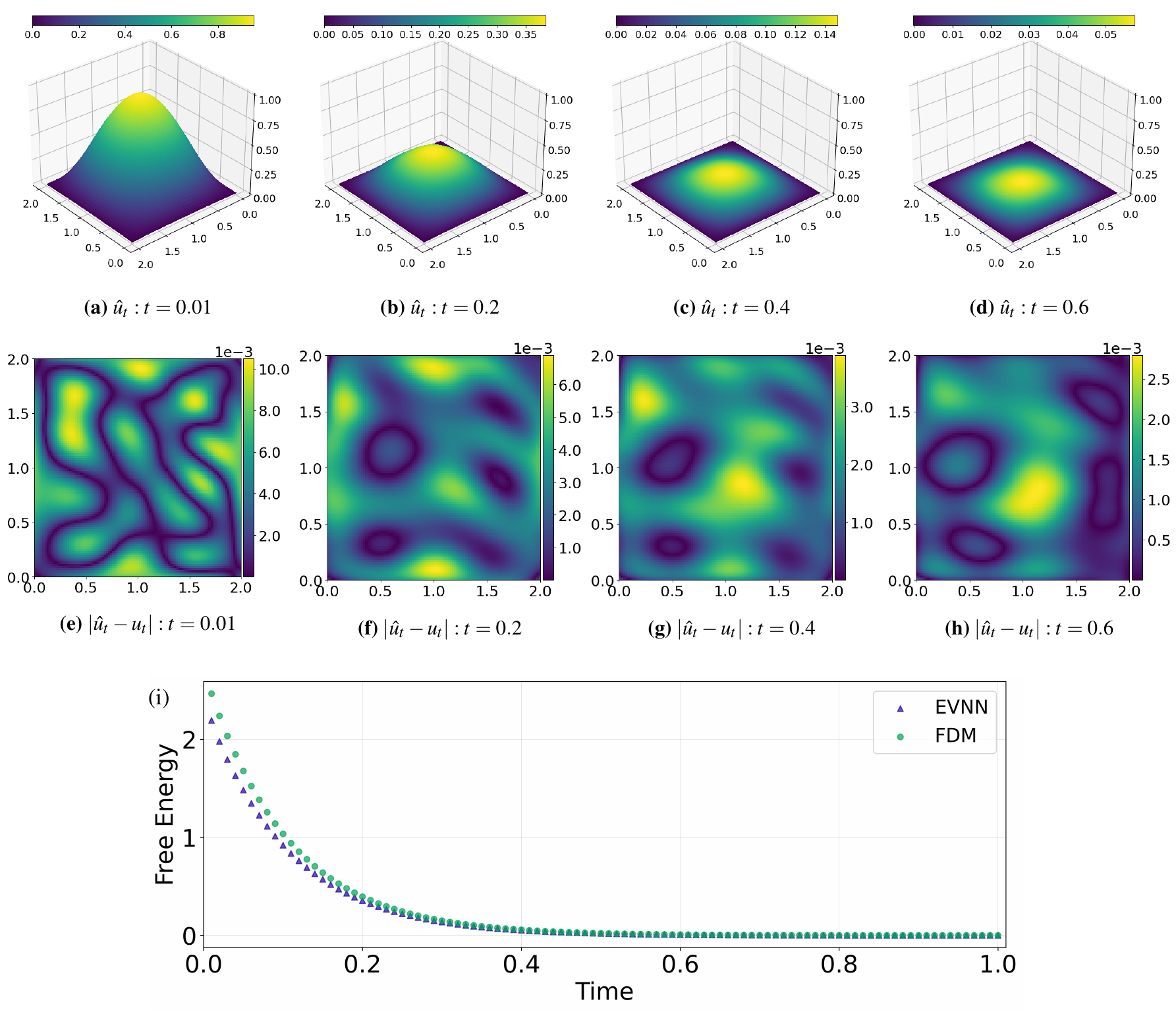}
\caption{Numerical results for the heat equation.  (a) - (d) are the neural network solutions at $t = 0.01, 0.2, 0.4$ and $0.6$. (e) - (h) are the absolute differences between the EVNN solutions and the FDM solutions. (i) The evolution of discrete free energy with respect to time for both the EVNN and FDM solutions.}
\label{fig:heat}
\end{figure}

In the numerical simulation, we take $\tau = 0.01$ and use a 1-block ResNet with $20$ nodes in each layer. The nonlinear activation function is chosen to be $\tanh$. The total number of parameters is 501. To achieve energy stability, we fix the training samples to eliminate the numerical fluctuation arising from estimating the integration in \eqref{eq:innerloop}. The set of training samples comprises a $301 \times 301$ uniform grid on $(0, 2)^2$, and an additional $1000$ uniformly spaced grid points on each edge of the boundary. To test the accuracy of the solution over time, we compare the neural network solution with a numerical solution obtained by a finite difference method (FDM). We apply the standard central finite difference and an implicit Euler method to obtain the FDM solution. The numerical scheme can be reformulated as an optimization problem: \begin{equation}\label{FDM_op}
\phi_h^{n+1} = \mathop{\arg\min}_{u_h \in \mathcal{C}} \frac{1}{2 \tau} \| u_h - u_h^n \|^2_h +  \mathcal{F}_h (u_h), \quad \mathcal{F}_h (u_h) = \frac{1}{2} \langle  \nabla_h u_h, \nabla_h u_h  \rangle_{*}~.
\end{equation}
Here, $\mathcal{C} = \{  u_{i, j} ~|~  0 < i < I, \quad 0 < j < J \}$ is the set of grid functions defined on a uniform mesh, $\nabla_h$ is the discrete gradient, $\mathcal{F}_h$ is the discrete free energy, $\| \cdot \|_h$ is the discrete $L_2$-norm induced by the discrete $L_2$-inner product $\langle f, g \rangle = h^2 \sum_{i=1}^{N-1} f_{i,j} g_{i,j}$, and $\langle \cdot, \cdot \rangle_{*}$ is the discrete inner product defined on the staggered mesh points \cite{liu2021structure, salgado2022classical}.
We use a $101 \times 101$ uniform grid and take $\tau = 0.01$ in the FDM simulation. Figure~\ref{fig:heat}(a)-(d) shows the NN solution at $t = 0.01, 0.2, 0.4$ and $0.6$. The corresponding point-wise absolute errors are shown in Fig.~\ref{fig:heat}(e)-(h). 
The evolution of discrete free energy (evaluated on $101 \times 101$ uniform grid) for both methods is shown in Fig.~\ref{fig:heat}(i). Clearly, the neural-network-based algorithm achieves energy stability and the result is consistent with the FDM. It is worth mentioning that the size of the optimization problem in the neural network-based algorithm is much smaller than that in the FDM, although (\ref{FDM_op})
is a quadratic optimization problem and can be solved directly by solving a linear equation of $u_{i, j}$.

\subsubsection{Allen-Cahn equation}
Next, we consider an Allen-Cahn type equation, which is extensively utilized in phase-field modeling, and becomes a versatile technique to solve interface problems arising from different disciplines \cite{du2020phase}.
In particular, we focus on the following Allen-Cahn equation on $\Omega = (-1, 1)^2$:
\begin{equation}
\left\{
\begin{aligned}
& \tfrac{\partial \varphi}{\partial t}(\x,t)  = - \left( \tfrac{F'(\varphi)}{\varepsilon^2} - \triangle \varphi(\x,t) + 2 W \left( \textstyle \int \varphi \dd \x - A\right) \right)~,~~~~\x\in \Omega~, ~~~~ t >0~, \\
% & \mathbf{n}\cdot \nabla \phi(\x,t) = 0~,~~~~ \x \in \partial \Omega~.
& \varphi(\x, t) = -1~, \quad \x \in \pp \Omega ~, ~~~~ t >0~,\\
& \varphi(\x, 0) = - \tanh (10 ( \sqrt{4 x_1^2 + x_2^2} - 0.5 ) ) ~,~~~~\x\in \Omega~.\\
\end{aligned}
\right.
\label{eq:Allen-Cahn}
\end{equation}
where $\varphi(\x, t)$ is the phase-field variable, constants $\epsilon$, $W$ and $A$ are model parameters. 
%v The corresponding energy-dissipation law is given by
The Allen--Cahn equation (\ref{eq:Allen-Cahn}) can be regarded as an $L^2$-gradient flow,  derived from an energy-dissipation law
\begin{equation}\label{AC_Energy}
\textstyle \frac{\dd}{\dd t} \mathcal{F}[\varphi] = - \int_{\Omega} |\varphi_t|^2 \dd \x, \quad  \mathcal{F}[\varphi] = \int_{\Omega }\frac{1}{2}\left|\nabla\varphi\right|^2 + \frac{1}{4\varepsilon^2}\left(\varphi^2 - 1\right)^2 \dd \x  + W\left( \textstyle \int_{\Omega} \varphi \dd \x - A\right)^2.
\end{equation}
Here, the last term in the energy is a penalty term for the volume constraint $\int \varphi \dd \x = A$, as the standard Allen-Cahn equation does not preserve the volume fraction $\int \varphi \dd \x$ over time. In the numerical simulation, we take $A = - (4 - \pi r^2) + \pi r^2, r = 0.5,  \frac{1}{\epsilon^2} = 100$, and $W = 1000$.

\begin{figure}[!htb]
  \centering
    % \caption*{$t = 0$}
    \includegraphics[width = 0.95 \linewidth]{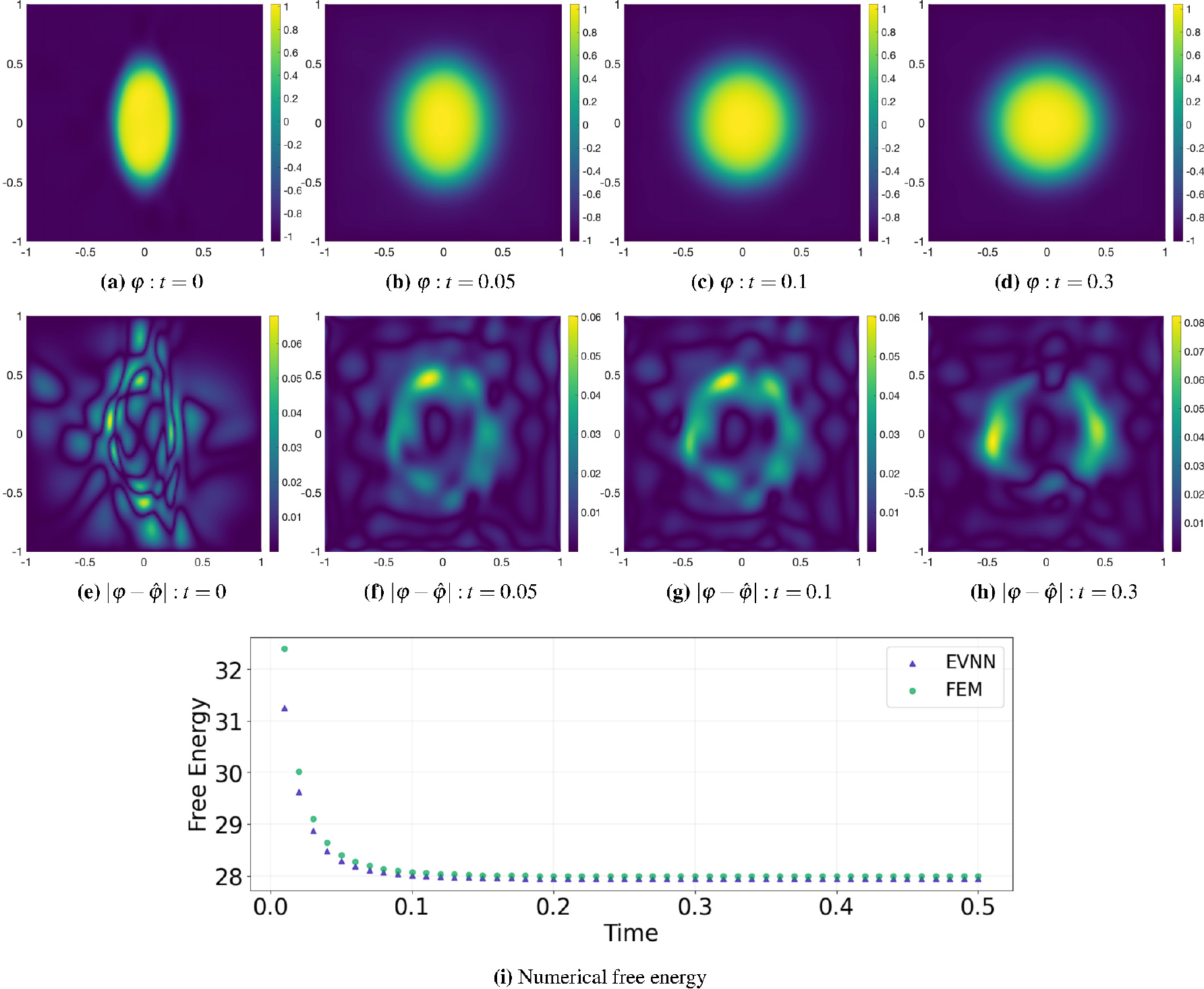}
 \caption{Numerical results for the Allen--Cahn equation with a volume constraint. (a)-(d) NN solutions at $t = 0, 0.05, 0.1$ and $0.3$, respectively. (e)-(h) Absolute differences between the EVNN solutions and the FEM solutions at $t = 0, 0.05, 0.1$, and $0.3$, respectively. (i) Evolution of numerical free energies with respect to time for both the FEM and EVNN solutions.}
  \label{fig:PhaseField}
\end{figure}

Due to the complexity of the initial condition, we utilized a larger neural network compared to the one employed in the previous subsection. Our choice was a 1-block ResNet, which has two fully connected layers, each containing 20 nodes. The total number of parameters is $921$. The nonlinear activation function is chosen to be $\tanh$. The set of training samples comprises a $301 \times 301$ uniform grid on $(-1, 1)^2$, and an additional $1000$ uniformly spaced grid points on each edge of the boundary. To test the numerical accuracy of the EVNN scheme for this problem, we also solve Eq.~(\ref{eq:Allen-Cahn}) by a finite element method, which approximates the phase-field variable $\varphi$ by a piece-wise linear function $\varphi_h(\x, t) = \sum_{i=1}^N \gamma_i(t) \psi_i(\x)~,$ where $\psi_i(X)$ are hat functions supported on the computational mesh. Inserting $\varphi_h$ into the continuous energy--dissipation law (\ref{AC_Energy}), we get a discrete energy--dissipation law with the discrete energy and dissipation given by $$\mathcal{F}_h^{\rm{FEM}} ({\bm \gamma}) = \textstyle \sum_{e = 1}^{N_e} \int_{\tau_e} \frac{1}{2} \Big| \sum_{i=1}^N \gamma_i \nabla \psi_i(\x) \Big|^2 + \frac{1}{\epsilon^2} \sum_{i=1}^N (\gamma_i^2 - 1)^2 \psi_i(\x)  \dd \x~,$$
and $\mathcal{D}_h^{\rm{FEM}} ({\bm \gamma}, {\bm \gamma}') = \sum_{e = 1}^{N_e} \int_{\tau_e} \Big| \sum_{i=1}^N \gamma_i'(t) \psi_i(\x) \Big|^2 \dd \x$
respectively. Here $\tau_e$ is used to denote a finite element cell, and $N_e$ is the number of cells. %We also construct a piece-wise linear approximation to the nonlinear term in the discrete energy. 
This form of discretization was used in \cite{xu2019stability}, which constructs a piece-wise linear approximation to the nonlinear term in the discrete energy.
We can update $\gamma_i$ at each time step by solving the optimization problem (\ref{Min_move_theta}), i.e.,
\begin{equation}\label{op_AC}
{\bm \gamma}^{n+1} = \mathop{\arg\min}_{{\bm \gamma} \in \mathbb{R}^N} \frac{1}{2 \tau}   {\sf D} ({\bm \gamma} - {\bm \gamma}^n) \cdot ({\bm \gamma} - {\bm \gamma}^n) + \mathcal{F}_h ({\bm \gamma})\ , 
\end{equation}
where ${\sf D}_{ij} = \int_{\Omega} \psi_i \psi_j \dd \x $ is the mass matrix. The optimization problem (\ref{op_AC}) is solved by L-BFGS in our numerical implementation.

The simulation results are summarized in Fig.~\ref{fig:PhaseField}. It clearly shows that our method can achieve comparable results with the FEM.  Numerical simulation of phase-field type models is often challenging. To capture the dynamics of the evolution of the thin diffuse interface, the mesh size should be much smaller than $\epsilon$, the width of the diffuse interface. Traditional numerical methods often use an adaptive or moving mesh approach to overcome this difficulty. In contrast, a neural network-based numerical scheme has a mesh-free feature. The number of parameters of the neural network can be much smaller than the number of samples needed to resolve the diffuse interface. Consequently, the dimension of the optimization problem in the NN-based scheme is much smaller than in the FEM scheme without using adaptive or moving meshes.

\subsection{Generalized diffusions}
In this subsection, we apply the proposed Lagrangian EVNN scheme to solving a Fokker-Planck equation and a porous medium equation.
\subsubsection{Fokker-Planck equation}
We first consider a Fokker-Planck equation
\begin{equation}
\label{eq:fokkersol0}
\begin{cases}
    & \rho_t =  \nabla \cdot ( \nabla \rho + \rho \nabla V)~, \quad \x  \in  \mathbb{R}^d~, \quad t\in (0, T]~,  \\
    & \rho(\x, 0) = \rho_0(\x)~, \quad \quad \quad \quad \x \in \mathbb{R}^d~, \\
\end{cases}
\end{equation}
where $V(\x)$ is a prescribed potential energy. The Fokker-Planck equation can be viewed as a generalized diffusion satisfying the energy-dissipation law
\begin{equation}
    \label{eq:fokkersol}
    \frac{\dd}{\dd t} \int_{\mathbb{R}^d} \rho \ln \rho + \rho V(\x) \dd \x = - \int_{\mathbb{R}^d} \rho |\uvec|^2 \dd \x~,
\end{equation}
where the probability density $\rho$ satisfies the kinematics $\rho_t + \nabla \cdot (\rho \uvec) = 0$.  

We test our numerical scheme for solving the Fokker-Planck equation in 2D and 4D, respectively. For both numerical experiments, we adopt the augmented-ICNN proposed in~\cite{huang2020convex}. We use a $1$ ICNN block of $6$ fully connected layers with $32$ hidden nodes in each layer. The activation function is chosen as the Gaussian-Softplus function $\sigma(x) = \textstyle \sqrt{\tfrac{2}{\pi}}(\sqrt{2}x \int_{0}^{x/\sqrt{2}} e^{-t^{2}} d t + e^{-\frac{x^2}{2}} + \sqrt{\tfrac{\pi}{2}}x )$. In addition, we use an L-BFGS to solve the optimization problem at each time step. 

\begin{figure}[tbh]
\includegraphics[width = 0.95 \linewidth]{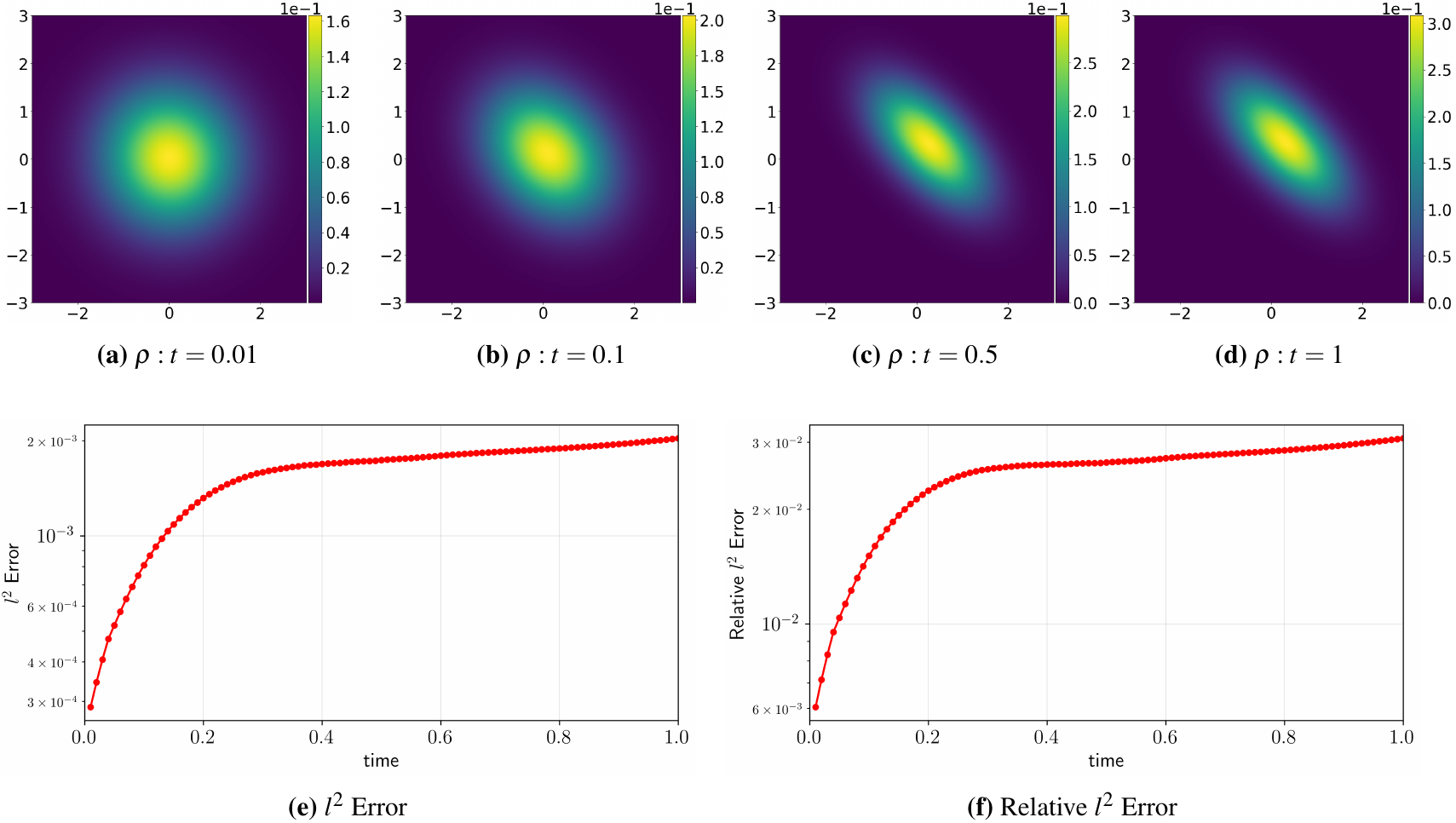}
\caption{Numerical results for the 2D Fokker--Planck equation. (a)-(d) Predicted solution on a $301 \times 301$ uniform mesh on $(-3, 3)^2$ at $t = 0, 0.1, 0.5$ and $1$, respectively. (e)-(f) The absolute and relative $l_2$-errors of the solution over time, evaluated on the uniform mesh.}
\label{fig:fokker}
\end{figure}

In the 2D case, we consider $V(\x)$ as follows:
\begin{equation*}
  \textstyle  V = \frac{1}{2}(\x - \mu_{target})^T\Sigma^{-1}_{    target} (\x - \mu_{target})
    ~~{\rm with}~~ \mu_{target} = \left(\frac{1}{3}, \frac{1}{3}\right)~,
\end{equation*}
and $\Sigma_{target} = \begin{bmatrix}
\frac{5}{8} & -\frac{3}{8} \\
-\frac{3}{8} & \frac{5}{8}
\end{bmatrix}$. The initial condition $\rho_0(\x)$ is set to be the 2D standard Gaussian $\mathcal{N}({\bf 0}, {\bf I})$. The exact solution of Eq.~(\ref{eq:fokkersol0}) takes the following analytical form \cite{hwang2021deep}
\begin{equation}
    \rho(\x, t) \sim \mathcal{N}(\mu(t), \Sigma(t))~,
\end{equation}
where $\mu(t) = (1 - \e^{-4t})\mu_{target}$ and $\Sigma(t) = \begin{bmatrix}
\frac{5}{8}+\frac{3}{8} \times e^{-8 t} & -\frac{3}{8}+\frac{3}{8} \times e^{-8 t} \\
-\frac{3}{8}+\frac{3}{8} \times e^{-8 t} & \frac{5}{8}+\frac{3}{8} \times e^{-8 t}
\end{bmatrix}$. We draw 10000 samples from $\rho_0(\x)$ as the training set. As an advantage of the neural network-based algorithm, we can compute $\rho^n(\x)$ point-wisely through a function composition $\rho^{n}(\x) = \rho_0 \tilde{\circ} (\Phi^{n})^{-1} (\x)$, where $(\Phi^{n})^{-1} (\x)$ can be computed by solving a convex optimization problem, i.e., $\Phi^{-1}(\x) = \mathop{\arg\min}_{{\bm y}} \Phi({\bm y}) - \x^{\rm T} {\bm y}$.
Fig.~\ref{fig:fokker}(a) - (d) shows the predicted density on a $301 \times 301$ uniform grid on $(-3, 3)^2$. 
The absolute and relative $l_2$-errors on the grid are shown in Figs.~\ref{fig:fokker}(e)-(f). It can be noticed that the relative $l^2$-error for the predicted solution is around $10^{-2}$, which is quite accurate given the limited number of samples used.

\begin{figure}[!t]
\centering
\includegraphics[width = 0.95 \linewidth]{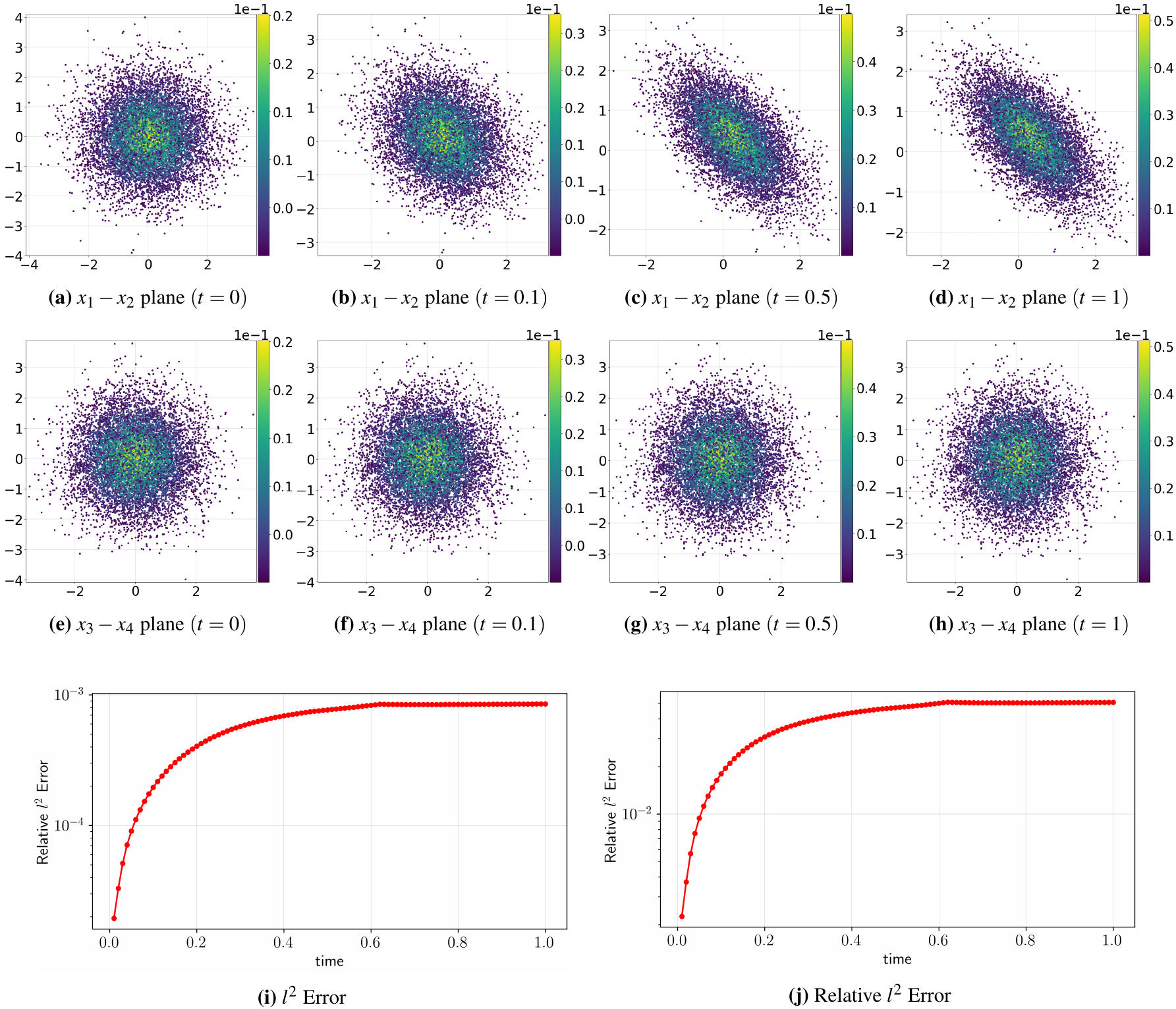}
\caption{Numerical results for the 4D Fokker--Planck equation. (a)-(d) 
Distribution of samples, projected in $x_1$-$x_2$ plane, as well as its weight at $t = 0, 0.1, 0.5$ and $1$, respectively. (e)-(h) Distribution of samples, projected in $x_3$-$x_4$ plane, as well as its weight at  $t = 0, 0.1, 0.5$ and $1$, respectively. (i)-(j) The $l_2$-errors and the relative $l^2$-errors of the solution over time, evaluated in the training set. }
\label{fig:fokker4d}
\end{figure}

In the 4D case, we take $V(\x)$ to be
\begin{equation*}
   V = \tfrac{1}{2}(\x - \mu_{target})^T\Sigma^{-1}_{target}(\x - \mu_{target})
    \quad {\rm with} \quad   \mu_{target} = \left(\tfrac{1}{3}, \tfrac{1}{3}, 0, 0\right)~,
\end{equation*}
and $\Sigma_{target} = \begin{bmatrix}
\frac{5}{8} & -\frac{3}{8} \\
-\frac{3}{8} & \frac{5}{8}
\end{bmatrix}\bigoplus\begin{bmatrix}
1 & 0 \\
0 & 1
\end{bmatrix}$. The initial condition $\rho_0(\x)$ is set to be a 4D standard normal distribution $\mathcal{N}({\bf 0}, {\bf I})$. The exact solution of Eq.~(\ref{eq:fokkersol0}) follows the normal distribution with 
$\mu(t) = \left((1 - \e^{-4t})\tfrac{1}{3}, (1 - \e^{-4t})\tfrac{1}{3}, 0, 0\right),$ and  $\Sigma(t) = \begin{bmatrix}
\frac{5}{8}+\frac{3}{8} \times e^{-8 t} & -\frac{3}{8}+\frac{3}{8} \times e^{-8 t} \\
-\frac{3}{8}+\frac{3}{8} \times e^{-8 t} & \frac{5}{8}+\frac{3}{8} \times e^{-8 t}
\end{bmatrix} \bigoplus \begin{bmatrix}
1 & 0 \\
0 & 1
\end{bmatrix}.$
Here $\bigoplus$ is the direct sum.
Same to the 2D case, we draw 10000 initial samples from $\rho_0(\x)$.
Fig.~\ref{fig:fokker4d} shows the evolution of samples as well the values of the numerical solution on individual samples over space and time. Interestingly, the relative $l^2$-error in 4D case is similar to the 2D case. The dimension-independent error bound suggests the potential of the current method for handling higher-dimensional problems.

Fokker--Planck type equations have a wide application in machine learning.
One of the fundamental tasks in modern statistics and machine learning is to
estimate or generate samples from a target distribution $\rho^*(\x)$, which might be completely known, partially known up to the normalizing constant, or empirically given by samples.
Examples include Bayesian inference \cite{blei2017variational},
numerical integration \cite{muller2019neural}, space-filling design \cite{pronzato2012design}, density estimation \cite{tabak2010density}, and generative learning \cite{papamakarios2021normalizing}. These problems can be transformed as an optimization problem, which is to seek for
%a set of particle $\{ \x_i \}_{i=1}^N$ 
a $\rho^{\rm opt} \in \mathcal{Q}$ by solving an optimization problem
\begin{equation}\label{op_pro}
\rho^{\rm opt} = \mathop{\arg\min}_{\rho \in \mathcal{Q}} D(\rho || \rho^*)~, % \quad \rho_N = \frac{1}{N} \sum_{i=1}^N \delta(\x - \x_i),
\end{equation}
where $\mathcal{Q}$ is the admissible set, $D(p || q)$ is a dissimilarity function that assesses the differences between two probability measures $p$ and $q$. The classical dissimilarities include the Kullback--Leibler (KL) divergence and the Maximum Mean Discrepancy (MMD). The optimal solution to the optimization problem can be obtained by solving a Fokker--Planck type equation, given by
\begin{equation}\label{Eq_G_FP}
    \rho_t = \nabla \cdot ( \rho \nabla \mu), \quad \mu = \frac{\delta D(\rho || \rho^*)}{\delta \rho}.
\end{equation}
The developed numerical approach has potential applications in these machine learning problems. 

\begin{figure}[!t]
  \centering
  \includegraphics[width = 0.95 \linewidth]{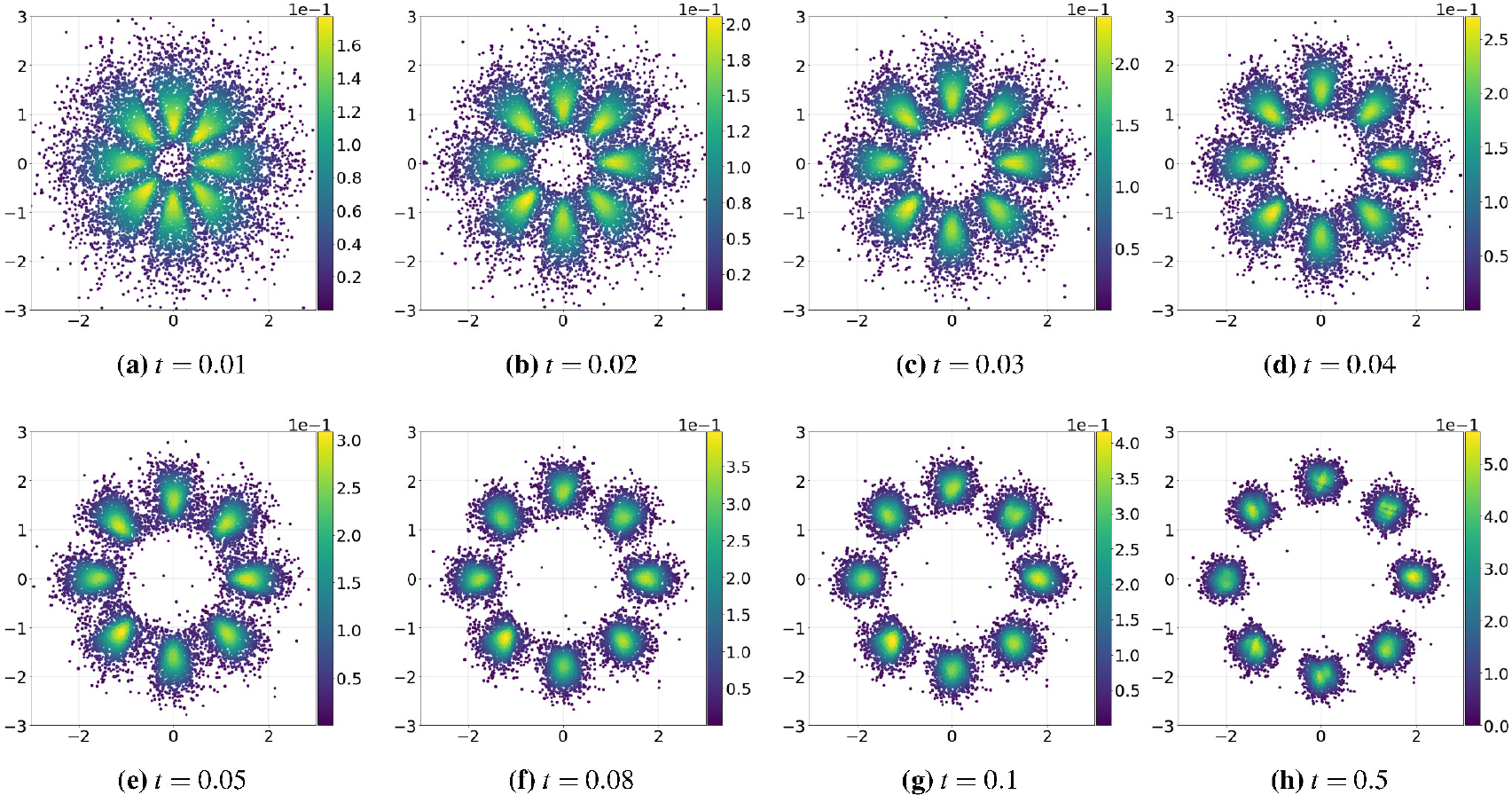}
  \caption{The distribution of samples and their weight at different times for the Fokker--Planck equation with a complicated $V(\x)$, corresponds to an eight-components mixture Gaussian distribution. The initial 10000 samples are drawn from a standard normal distribution $\mathcal{N}({\bf 0}, {\bf I})$. }
    % Numerical solutions for Fokker-Planck equation with complicated $V({\bm x)}$.}
  \label{fig:fokker8d}
  \end{figure}

To illustrate this point, we consider a toy problem that is widely used in the machine learning literature \cite{huang2020convex}. The goal is to sample from a target distribution $\rho^*$, which is known up to a normalizing constant.
We take the dissimilarity function as the KL divergence ${\rm KL}(\rho || \rho^*) =  \int_{\Omega}  \rho(\x) \ln \left( \frac{\rho}{\rho^*}  \right) \dd \x,$ then Eq. (\ref{Eq_G_FP}) is reduced to the Fokker--Planck equation (\ref{eq:fokkersol0}) with $V(\x) = - \ln \rho^*$. 
In the numerical experiment, we take the target distribution $\rho^{*}(\x)=\frac{1}{8}\sum_{i=1}^8 N(\bm x|\bm \mu_i,\bm \Sigma)$, an eight components mixture Gaussian distribution,
% \]
where $\bm \mu_1=(0,4)$, $\bm \mu_2=(2.8,2.8)$, $\bm \mu_3=(4,0)$, $\bm \mu_4=(-2.8,2.8)$, $\bm \mu_5=(-4,0)$, $\bm \mu_6=(-2.8,-2.8)$, $\bm \mu_7=(0,-4)$, $\bm \mu_8=(2.8,-2.8)$, and $\bm \Sigma = \mathrm{diag}(0.2, 0.2)$.  %\left[\begin{array}{rr} 0.2, & 0 \\ 0, & 0.2 \end{array}\right]$. 
The simulation results are summarized in Fig.~\ref{fig:fokker8d}, in which we first draw 10000 samples from a standard normal distribution $\mathcal{N}({\bf 0}, {\bf I})$ and show the distribution of samples as well as their weights at different time.  It can be noticed that the proposed neural network-based algorithm can generate a weighted sample for a complicated target distribution.

\subsubsection{Porous medium equation}
Next, we consider a porous medium equation (PME), $\rho_t = \Delta \rho^{\alpha}~,$where $\alpha > 1$ is a constant. The PME is a typical example of nonlinear diffusion equations.
One important feature of the PME is that the solution to the PME has a compact support at any time $t > 0$ if the initial data has a compact support. The free boundary of the compact support moves outward with a finite speed, known as the property of finite speed propagation \cite{vazquez2007porous}. As a consequence, numerical simulations of the PME are often difficult by using Eulerian methods, which may fail to capture the movement of the free boundary and suffer from numerical oscillations \cite{liu2020lagrangian}. In a recent work \cite{liu2020lagrangian}, the authors developed a variational Lagrangian scheme using a finite element method. Here we show the ability of the Lagrangian EVNN scheme to solve the PME with a free boundary.
Following \cite{liu2020lagrangian}, we employ the energy-dissipation law 
\begin{equation}
    \frac{\dd}{\dd t} \int_{\mathbb{R}^d} \frac{\alpha}{(\alpha - 1)(\alpha-2)} \rho^{\alpha-1} \dd \x = - \int_{\mathbb{R}^d} |\uvec|^2 \dd \x~
    \end{equation}
to develop the EVNN scheme.

\begin{figure}[!tbh]
\centering
\includegraphics[width = 0.95 \linewidth]{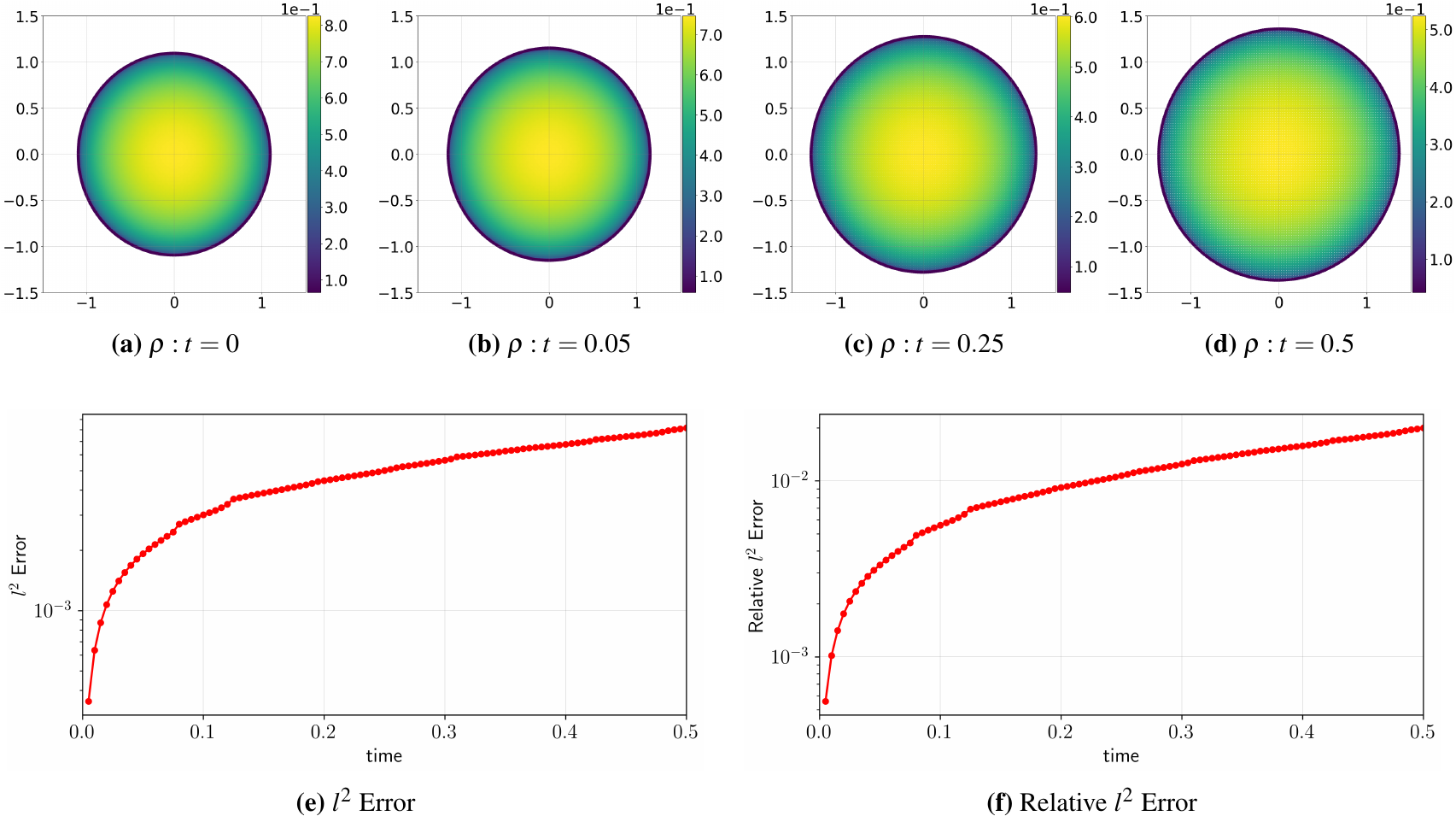}
\caption{Numerical results for the porous medium equation. (a)-(d) The numerical solution  at $t = 0, 0.05, 0.25$, and $0.5$, respectively. (e)-(f) The $L_2$-errors of the solution.}
\label{fig:porus}
\end{figure}

We take $\alpha = 4$ in the simulation.
To test the numerical accuracy of the EVNN scheme for solving the PME, we consider a 2D Barenblatt--Pattle solution to the PME. 
The Barenblatt--Pattle solution in $d$-dimensional space is given by
\begin{equation}\label{2D_B}
B_{\alpha} (\x ,t) = \textstyle t^{-k} \bigl[ \bigl(C_0 - \frac{k(\alpha-1)}{2 d \alpha} \frac{|\x|^2}{t^{2k/d}}\bigr)_{+}  \bigr]^{1/(\alpha-1)}~, \quad \x \in \mathbb{R}^d~,
\end{equation}
where $k = (\alpha - 1 + 2/d)^{-1}$, $u_{+} := \max(u, 0)$, and $C_0$ is a constant that is related to the initial mass. This solution is radially symmetric, self-similar, and has compact support  $|\x| \leq \xi_{\alpha}(t)$ for any finite time, where
$\xi_{\alpha}(t) = \textstyle \sqrt{\frac{2d \alpha C_0}{k (\alpha - 1)}} t^{k/d}$.
We take the Barenblatt solution (\ref{2D_B}) with $C_0 = 0.1$ at $t = 0.1$ as the initial condition, and compare the numerical solution at $T$ with the exact solution $B_{\alpha}(\x, T + 0.1)$. 
We choose $\{ x_i^0 \}_{i=1}^N$ as a set of uniform grid points 
inside a disk $\{ (x, y) ~|~ \sqrt{x^2 + y^2} \leq \xi_{\alpha}(0.1) \}$. To visualize the numerical results, we also draw $500$ points distributed with
equal arc length on the initial free boundary and evolve them using the neural network.
The numerical results with $\tau = 0.005$ are shown in Fig.~\ref{fig:porus}. 
It can be noticed that the proposed neural network-based algorithm can well approximate the Barenblatt--Pattle solution and capture the movement of the free boundary.

We note that the numerical methods to solve generalized diffusions (\ref{ED_Diff}) can also be formulated in the Eulerian frame of reference based on the notion of Wasserstein gradient flow \cite{jordan1998variational}. % We refer the interested readers to a recent work  \cite{hwang2021deep} for a neural network implementation of the minimizing movement scheme with the Wasserstein distance.
However, it is often challenging to compute the Wasserstein distance between two probability densities efficiently with high accuracy, which often requires solving an additional min-max problem or minimization problem \cite{benamou2000computational, benamou2016augmented, carrillo2022primal, hwang2021deep}. In a recent study \cite{hwang2021deep}, the authors initially employ a fully connected neural network to approximate $\rho$ and subsequently utilize an additional ICNN to calculate the Wasserstein distance between two probability densities. 
Compared with their approach, our method is indeed more efficient and accurate.

\section{Conclusion}
\label{sec:con}

In this paper, we develop structure-preserving EVNN schemes for simulating the $L^2$-gradient flows and the generalized diffusions by utilizing a neural network as a tool for spatial discretization, within the framework of the discrete energetic variational approach. These numerical schemes are directly constructed 
based on a prescribed continuous energy-dissipation law for the system, without the need for the underlying PDE. %
The incorporation of mesh-free neural network discretization opens up exciting possibilities for tackling high-dimensional gradient flows arising in different applications in the future.

Various numerical experiments are presented to demonstrate the accuracy and energy stability of the proposed numerical scheme. In our future work, we will explore the effects of different neural network architectures, sampling strategies, and optimization methods, followed by a detailed numerical analysis. Additionally, we intend to employ the EVNN scheme to investigate other complex fluid models, including the Cahn--Hilliard equation and Cahn--Hilliard--Navier--Stokes equations, as well as solve machine learning problems such as generative modeling and density estimation.

\section*{Acknowledgments}
C. L and Y. W were partially supported by NSF DMS-1950868 and DMS-2153029. Z. X was partially supported by
the NSF CDS\&E-MSS 1854779.

\noindent 

\bibliography{references}

% \bibliographystyle{asa}
% \bibliography{reference}

\end{document}